\documentclass[10pt]{article}
\usepackage[square,authoryear]{natbib}
\usepackage{marsden_article}
\usepackage{epstopdf}
\usepackage{setspace}
\usepackage{amscd}

\usepackage{caption2}

\DeclareGraphicsExtensions{.eps,.ps,.pdf}

\begin{document}

 \newtheorem{thm}{Theorem}[section]
 \newtheorem{cor}[thm]{Corollary}
 \newtheorem{lem}[thm]{Lemma}
 \newtheorem{prop}[thm]{Proposition}
 \newtheorem{defn}[thm]{Definition}
 \newtheorem{rem}[thm]{Remark}
 \numberwithin{equation}{section}

%
%
%
%
%
%
%

\title{\Large{\textbf{COVARIANT BALANCE LAWS IN CONTINUA WITH
MICROSTRUCTURE}}}

\author{Arash Yavari\thanks{School of Civil and Environmental Engineering,
  Georgia Institute of Technology, Atlanta, GA 30332. E-mail: arash.yavari@ce.gatech.edu. Research supported by the Georgia Institute of Technology.}
  \and Jerrold E. Marsden\thanks{Control and Dynamical Systems,
  California Institute of Technology, Pasadena, CA 91125. Research partially supported by the California Institute of Technology and NSF-ITR Grant ACI-0204932.} }

\maketitle

\begin{abstract} The purpose of this paper is to extend the Green-Naghdi-Rivlin
balance of energy method to continua with microstructure. The key
idea is to replace the group of Galilean transformations with the
group of diffeomorphisms of the ambient space. A key advantage is
that one obtains in a natural way all the needed balance laws on
both the macro and micro levels along with two Doyle-Erickson
formulas.

We model a structured continuum as a triplet of Riemannian
manifolds: a material manifold, the ambient space manifold of
material particles and a director field manifold. The
Green-Naghdi-Rivlin theorem and its extensions for structured
continua are critically reviewed. We show that when the ambient
space is Euclidean and when the microstructure manifold is the
tangent space of the ambient space manifold, postulating a single
balance of energy law and its invariance  under time-dependent
isometries of the ambient space, one obtains conservation of mass,
balances of linear and angular momenta but {\it not} a separate
balance of linear momentum.

We develop a covariant elasticity theory for structured continua
by postulating that energy balance is invariant under
time-dependent spatial diffeomorphisms of the ambient space, which
in this case, is the product of two Riemannian manifolds. We then
introduce two types of constrained continua in which
microstructure manifold is linked to the reference and ambient
space manifolds. In the case when at every material point, the
microstructure manifold is the tangent space of the ambient space
manifold at the image of the material point, we show that the
assumption of covariance leads to balances of linear and angular
momenta with contributions from both forces and micro-forces along
with two Doyle-Ericksen formulas. We show that   generalized
covariance  leads to two balances of linear momentum and a single
coupled balance of angular momentum.

Using this theory, we covariantly obtain the balance laws for two
specific examples, namely elastic solids with distributed voids
and mixtures. Finally, the Lagrangian field theory of structured
elasticity is revisited and a connection is made between
covariance and Noether's theorem.
\end{abstract}

\begin{description}
\item[Keywords:] Continuum Mechanics, Elasticity, Generalized
Continua, Couple Stress, Energy Balance.
\end{description}

\tableofcontents


\vskip 0.4 in
\section{Introduction}

The idea of generalized continua goes back to the work of
\cite{Cosserat1909}. The main idea in generalized continua is to
consider extra degrees of freedom for material points in order to
be able to better model materials with microstructure in the
framework of continuum mechanics. Many developments have been
reported since the seminal work of the Cosserat brothers.
Depending on the specific choice of kinematics, generalized
continua are called polar, micropolar, micromorphic, Cosserat,
multipolar, oriented, complex, etc. (see \cite{Green64a},
\cite{KafEr1971}, \cite{Toupin1962}, \cite{Toupin1964},
\cite{Mindlin1964} and references therein). The more recent
developments can be seen in \cite{Car1989}, \cite{CarMar2003},
\cite{FabMar2005}, \cite{EpsteindeLeon1998}, \cite{Muschik2001},
\cite{Slawianowski2005} and references therein. For a recent
review see \cite{MarianoStazi2005}.

By choosing a specific form for the kinetic energy density of
directors, \citet{Cov1975} obtained the balance laws of a Cosserat
continuum with three directors by imposing invariance of energy
balance under rigid translations and rotations in the current
configuration. A similar work was done by \citet{Buggisch1973}.
\citet{CarPodWil1982} obtained the balance laws for a continuum
with the so-called affine microstructure by postulating invariance
of balance of energy under time-dependent rigid translations and
rotations of the deformed configuration. The main assumption there
is that the orthogonal second-order tensor representing the affine
microdeformations remains unchanged under a rigid translation but
is transformed liked a two-point tensor under a rigid rotation in
the deformed configuration. Accepting this assumption, one obtains
conservation of mass, the standard balance of linear momentum and
balance of angular momentum, which in this case states that the
sum of Cauchy stress and some new terms is symmetric. Recently,
\citet{FabMar2005} conducted an interesting study of the geometric
structure of complex continua and studied different geometric
aspects of continua with microstructure. \citet{CarMar2003}
studied the Lagrangian field theory of Coserrat continua and
obtained the Euler-Lagrange equations for standard and
microstructure deformation mappings. However, in their Lagrangian
density they did not consider an explicit dependence on the metric
of the order-parameter manifold. In this paper, we will consider
an explicit dependence of the Lagrangian density on metrics of
both standard and microstructure manifolds. One should remember
that the original developments in the theory of generalized
continua in the Sixties were variational \citep{Toupin1962,
Toupin1964}. However, revisiting the Lagrangian field theory of
structured continua in the language of modern geometric mechanics
may be worthwhile.

It is believed that kinematics of a structured continuum can be
described by two independent maps, one mapping material points to
their current positions and one mapping the material points to
their directors \citep{MaHu1983}. Looking at the literature one
can see that for a Cosserat continuum (and even for multipolar
continua \citep{Green64, Green64a}), the only balance laws are the
standard balances of linear and angular momenta; couple stresses
do not enter into balance of linear momentum but do enter into
balance of angular momentum and make the Cauchy stress
unsymmetric. This is indeed different from the situation in the
so-called complex continua or continua with microstructure
\citep{Car1989, CarMar2003, FabMar2005}, where one sees separate
balance laws for microstresses. \cite{MaHu1983} postulated two
balances of linear momenta. However, it is not clear why, in
general, one should see two balances of linear momentum and only
one balance of angular momentum. In other words, why do standard
and microstructure forces interact only in the balance of angular
momentum? It should be noted that in all the existing
generalizations of Green-Naghdi-Rivlin (GNR) Theorem (see
\cite{Green64}) to generalized continua the standard Galilei group
$\mathcal{G}$ is considered. It is always assumed that rigid
translations leave the micro-kimenatical variables and their
corresponding forces unchanged (with no rigorous justification)
and these quantities come into play only when rigid rotations are
considered.

It is known that the traditional formulation of balance laws of
continuum mechanics are not intrinsically meaningful and heavily
depend on the linear structure of Euclidean space.
\citet{MaHu1983} resolved this shortcoming of the traditional
formulation by postulating a balance of energy, which is
intrinsically defined even on manifolds, and its invariance under
spatial changes of frame. This results in conservation of mass,
balance of linear and angular momenta and the Doyle-Ericksen
formula. Similar ideas had been proposed in \cite{Green64} for
deriving balance laws by postulating energy balance invariance
under Galilean transformations. For more details and discussions
on material changes of frame see \citet{YaMaOr2006}. See also
\cite{Ya2008}, \cite{YavariOzakin2008}, and \cite{YaMa2008} for
similar discussions. A natural question to ask is whether it is
possible to develop covariant theories of elasticity for
structured continua. As we will see shortly, the answer is
affirmative.

Similar to Noether's theorem that makes a connection between
conserved quantities and symmetries of a Lagrangian density, GNR
theorem makes a connection between balance laws and invariance
properties of balance of energy. One major difference between the
two theorems is that in GNR theorem one looks at balance of energy
for a finite subbody, i.e., a global quantity, and its invariance,
while in Noether's theorem symmetries are local properties of the
Lagrangian density.

In some applications, e.g., recent applications of continuum
mechanics to biology, one may need to enlarge the configuration
manifold of the continuum to take into account the fact that
changes in material points, e.g., rearrangements of
microstructure, etc., should somehow be considered in the
continuum theory, at least in an average sense. This was a
motivation for various developments for generalized continuum
theories in the last few decades. In a structured continuum, in
addition to the standard deformation mapping, one introduces some
extra fields that represent the underlying microstructure. In the
nondissipative case, assuming the existence of a Lagrangian
density that depends on all the fields, using Hamilton's principle
of least action one obtains new Euler-Lagrange equations
corresponding to microstructural fields \citep{Toupin1962,
Toupin1964, CarMar2003}. However, to our best knowledge, it is not
clear in the literature how one can obtain these extra balance
laws by postulating a single energy balance and its invariance
under some groups of transformations. This is the main motivation
of the present work.

To summarize, looking at the literature of generalized continua,
one sees that the structure of balance laws is not completely
clear. It is observed that there is always a standard balance of
linear momentum with only macro-quantities and a balance of
angular momentum, which has contributions from both macro- and
micro-forces. In some treatments there is no balance of
micro-linear momentum (see \cite{Toupin1962, Toupin1964,
CarPodWil1982, Ericksen1961}) while sometimes there is one, as in
\cite{GreenNaghdi95a, Car1989}. In particular, we can mention the
work of \citet{Leslie1968} on liquid crystals in which he starts
by postulating a balance of energy and a linear momentum balance
for micro-forces. In his work, he realizes that the balance of
micro-linear momentum cannot be obtained from invariance of energy
balance. To date, there have been several works on relating
balance laws of structured continua to invariance of energy
balance under some group of transformations. These efforts will be
reviewed in detail in the sequel.

This paper is organized as follows. In \S2 geometry of continua
with microstructure is discussed. \S3 discusses the previous
efforts in generalizing Green-Naghdi-Rivlin Theorem for
generalized continua. Assuming that the ambient space is Euclidean
and assuming that the microstructure manifold at every material
point is the tangent space of $\mathbb{R}^3$ at the spatial image
of the material point, we generalize GNR theorem. \S4 develops a
covariant theory of elasticity for those structured continua for
which microstructure manifold is completely independent of the
ambient space manifold in the sense that ambient space and
microstructure manifolds can have separate changes of frame. We
then develop a covariant theory of elasticity for those structured
continua in which microstructure manifold is somewhat linked to
the ambient space manifold. In particular, we study the case where
microstructure manifold is the tangent bundle of the ambient space
manifold. We also introduce a generalized notion of covariance in
which one postulates energy balance invariance under two
diffeomorphisms that act separately on micro and macro quantities
simultaneously. We study consequences of this generalized
covariance. In \S5, we look at two concrete examples of structured
continua, namely elastic solids with distributed voids and
mixtures. In both cases, we obtain the balance laws covariantly.
\S6 presents a Lagrangian field theory formulation of structured
continua. Noether's theorem and its connection with covariance is
also investigated. Concluding remarks are given in \S7.

\vskip 0.2 in
\section {\textbf{Geometry of Continua with Microstructure}}

A structured continuum is a generalization of a standard continuum
in which the internal structure of the material points is taken
into account by assigning to them some independent internal
variables or order parameters. For the sake of simplicity, let us
assume that each material point $\mathbf{X}$ has a corresponding
microstructure (director) field $\mathbf{p}$, which lies in a
Riemannian manifold $(\mathcal{M},\mathbf{g}_{\mathcal{M}})$. Note
that $\mathbf{p}$, in general, could be a tensor field. In
general, one may have a collection of director fields and the
microstructure manifold may not be Riemannian. However, these
assumptions are general enough to cover many problems of interest.
In this case our structured continuum has a configuration manifold
that consists of a pair of mappings
$(\varphi_t,\widetilde{\varphi}_t)$ \citep{MaHu1983, FabMar2005},
where $\mathbf{x}=\varphi_t(\mathbf{X})$ represents the standard
motion and $\mathbf{p}=\widetilde{\varphi}_t(\mathbf{X})$ is the
motion of the microstructure. Both $\varphi_t$ and
$\widetilde{\varphi}_t$ are understood as fields.
\begin{figure}[hbt]
\vspace*{0.3in}
\begin{center}
\includegraphics[scale=0.7,angle=0]{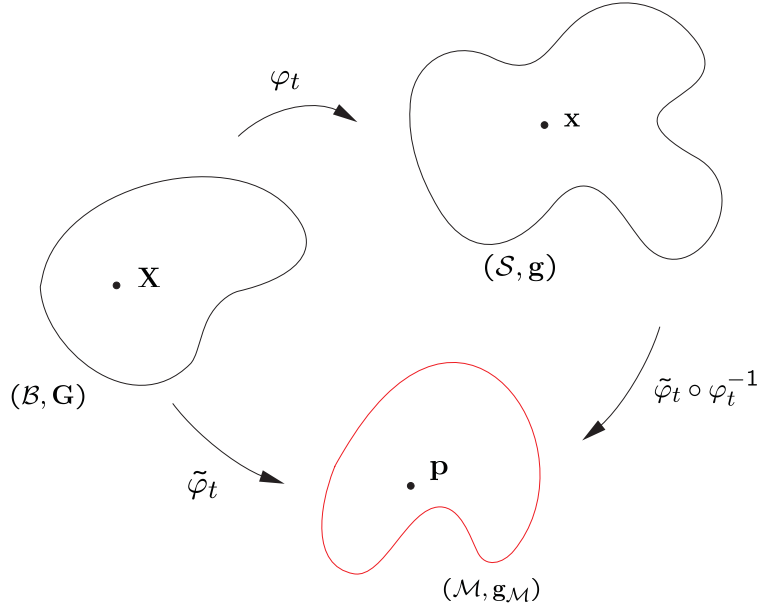}
\end{center}
\caption{\footnotesize Deformation mappings of a continuum with
microstructure.} \label{Cosserat_Geometry}
\end{figure}
As in the geometric treatment of standard continua, the current
configuration lies in an embedding space $\mathcal{S}$, which is a
Riemannian manifold with a metric $\mathbf{g}$. Note that ambient
space for the structured continuum is
$\overline{\mathcal{S}}=\mathcal{S}\times\mathcal{M}$ and for
every $\mathbf{X}\in\mathcal{B}$,
$\widetilde{\varphi}(\mathbf{X})$ lies in a separate copy of
$\mathcal{M}$. Here, we have assumed that the structured continuum
is microstructurally homogeneous in the sense that directors of
two material points $\mathbf{X}_1$ and $\mathbf{X}_2$ lie in two
copies of the same microstructure manifold $\mathcal{M}$ (see Fig.
\ref{Cosserat_Geometry}).
\begin{figure}[hbt]
\vspace*{0.3in}
\begin{center}
\includegraphics[scale=0.85,angle=0]{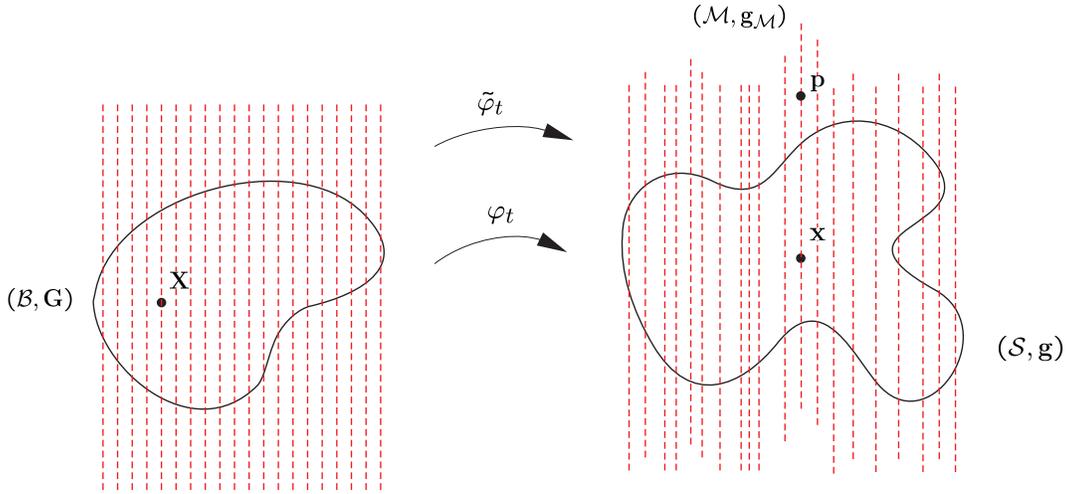}
\end{center}
\caption{\footnotesize Deformation of a continuum with
microstructure can be understood as a bundle map between two
trivial bundles. Here all is needed is the zero-section of the
reference bundle, i.e. the material manifold.}
\label{Cosserat_Bundle}
\end{figure}

More precisely, kinematics of a structured continuum is described
using fiber bundles (see, for instance, \cite{EpsteindeLeon1998}). Deformation of a
structured continua is a bundle map from the zero section of the
trivial bundle $\mathcal{B}\times\mathcal{M}_0$ (for some manifold
$\mathcal{M}_0$) to the trivial bundle
$\mathcal{S}\times\mathcal{M}$ (see Fig. \ref{Cosserat_Bundle}).
Corresponding to the two maps $\varphi_t$ and
$\widetilde{\varphi}_t$, there are two velocities, which have the
following material forms
\begin{equation}
   \mathbf{V}(\mathbf{X},t)=\frac{\partial \varphi_t(\mathbf{X})}{\partial t}\in T_{\mathbf{x}}\mathcal{S},
   ~~~\widetilde{\mathbf{V}}(\mathbf{X},t)=\frac{\partial \widetilde{\varphi}_t(\mathbf{X})}{\partial t}\in T_{\mathbf{p}}\mathcal{M}.
\end{equation}
Let us choose local coordinates $\{X^A\},\{x^a\},$ and
$\{p^{\alpha}\}$ on $\mathcal{B}$, $\mathcal{S}$ and
$\mathcal{M}$, respectively. In these coordinates
\begin{equation}
   \mathbf{V}(\mathbf{X},t)=V^a\mathbf{e}_a,
   ~~~\widetilde{\mathbf{V}}(\mathbf{X},t)=\widetilde{V}^{\alpha}~\widetilde{\mathbf{e}}_{\alpha},
\end{equation}
where $\{\mathbf{e}_a\}$ and $\{\widetilde{\mathbf{e}}_{\alpha}\}$
are bases for $T_{\mathbf{x}}\mathcal{S}$ and
$T_{\mathbf{p}}\mathcal{M}$, respectively, and
\begin{equation}
   V^a=\frac{\partial \varphi^a}{\partial t} ,~~~~~\widetilde{V}^{\alpha}=\frac{\partial \widetilde{\varphi}^{\alpha}}{\partial t}.
\end{equation}
In spatial coordinates
\begin{equation}
   \mathbf{v}(\mathbf{x},t)=\mathbf{V}\circ\varphi_t^{-1},~~~~~\widetilde{\mathbf{v}}(\mathbf{x},t)=\widetilde{\mathbf{V}}\circ\varphi_t^{-1}.
\end{equation}
In a local coordinate chart
\begin{equation}
   \mathbf{v}(\mathbf{x},t)=v^a\mathbf{e}_a,
   ~~~\widetilde{\mathbf{v}}(\mathbf{x},t)=\widetilde{v}^{\alpha}~\widetilde{\mathbf{e}}_{\alpha}.
\end{equation}
Here, for the sake of simplicity, we have assumed that our
structured continuum has one director field, which is assumed to
be a vector field. As was mentioned earlier, this is not the most
general possibility and in general one may need to work with
several director fields or even with a tensor-valued director
field. Generalization to these cases is straightforward.

\cite{MaHu1983} chose the classical viewpoint in taking
$\mathbb{R}^3$ to be the ambient space for material particles and
postulated the integral form of balances of linear and angular
momenta. 
The more natural approach would be to start from balance of energy
and look at consequences of its invariance under some
transformations. This is the approach we choose in this paper.
Note that the two maps $\varphi_t$ and $\widetilde{\varphi}_t$, in
general, are independent and interact only in the balance of
energy, i.e. power has contributions from both deformation maps.
The other important observation is that balance of energy is
written on an arbitrary subset
$\varphi_t(\mathcal{U})\subset\mathcal{S}$.

\vskip 0.2 in
\section {\textbf{The Green-Naghdi-Rivlin Theorem for a Continuum with Microstructure}}

In most theories of generalized continua, macro and micro-forces
enter the same balance of angular momentum because the ambient
space manifold and the manifold of microstructure are somewhat
related. Now the important question is the following: how can one
obtain two sets of balance of linear momentum, one for
micro-forces and one for marco-forces in such cases starting from
first principles? Of course, one can always postulate as many
balance laws as one needs in a theory. However, a fundamental
understanding of balance laws is crucial in any theory. Accepting
a Lagrangian viewpoint, one has two sets of Euler-Lagrange
equations as there are two independent macro and micro kinematic
variables (see \cite{Toupin1962, Toupin1964, FabMar2005}). Then,
assuming that these equations are satisfied, Noether's theorem
leads us to expect that any conserved quantity of the system
corresponds to some symmetry of the Lagrangian density. The
Lagrangian density can be invariant under groups of
transformations that act on the ambient and microstructure
manifolds simultaneously. For example, if one assumes that an
arbitrary element of $SO(3)$ acts simultaneously on $\mathcal{S}$
and $\mathcal{M}$ and Lagrangian density remains invariant, then
the conserved quantity is nothing but angular momentum with some
extra terms representing the effect of microstructure. However,
another possibility would be a symmetry in which an arbitrary
element of $SO(3)$ acts only on $\mathcal{M}$. Now one may ask why
the Lagrangian density should be invariant under simultaneous
actions of $SO(3)$ on $\mathcal{S}$ and $\mathcal{M}$.

A way out of this
difficulty may be to look for a generalization of the
Green-Naghdi-Rivlin theorem for continua with microstructure.
There have been several attempts in the literature to generalize
this theorem. In all the existing generalizations, it is assumed
that in a Galilean transformation, micro-forces and
micro-displacements remain unchanged under a rigid translation
while under a rigid rotation both micro and macro quantities
transform. Postulating invariance of balance of energy under an
arbitrary element of the Galilean group and accepting this
assumption, one obtains conservation of mass, the standard balance
of linear momentum and balance of angular momentum with some extra
terms that represent the effect of microstruture. However, this
does not give a micro-linear momentum balance. So, it is seen that
the link between energy balance invariance and balance of
micro-linear momentum is missing.

It should be noted that in most of the treatments of continua with
microstructure, the microstructure manifold $\mathcal{M}$ may not
be completely independent of the ambient space manifold
$\mathcal{S}$ and this may be a key point in understanding the
structure of balance laws. From a geometric point of view this
means that spatial and microstructure changes of frame may not be
independent, in general.

There have been several attempts in the literature to obtain
balance laws of generalized continua by energy invariance
arguments. \cite{CarPodWil1982} start from balance of energy and
postulate its invariance under rigid translations and rotations of
the current configuration. They assume that microstructure
quantities (kinematic and kinetic) remain unchanged under rigid
translations while under rigid rotations micro-forces transform
exactly like their macro counterparts. This somehow implies that
the microstructure manifold is not independent of the standard
ambient space. Under a rigid translation, each microstructure
manifold (fiber) translates rigidly and hence micro-forces and
directors remain unchanged. Under a rigid rotation directors and
their corresponding micro-forces transform exactly like their
macro counterparts because rotating a representative volume
element its director goes through the same rotation. This
invariance postulate results in the standard conservation of mass
and balance of linear and angular momenta. Balance of linear
momentum has its standard form while balance of angular momentum
has contributions from both forces and micro-forces. However, this
invariance argument does not lead to a separate balance of
micro-linear momentum.

\citet{Gurtin92} introduce a fine structure for each material
point. They then postulate two balances of energy, one in the
macro scale and one in the fine scale. The fine structure is
characterized by the limit $\epsilon\rightarrow 0$ of some scale
parameter $\epsilon$. Postulating invariance of these two balance
laws under rigid translations and rotations they obtain two sets
of balance of linear and angular momenta. They emphasize that
balance of micro-angular momentum only introduces a micro-couple
and offers nothing essential.

\cite{GreenNaghdi95a} and \citet{GreenNaghdi95b} start from
balance of energy and assume that it is invariant under the
transformation $\mathbf{v}\rightarrow\mathbf{v}+\mathbf{c}$, where
$\mathbf{v}$ is the spatial velocity field and $\mathbf{c}$ is an
arbitrary constant vector field. This gives the conservation of
mass and balance of linear momentum. Then they obtain a local form
for balance of energy and assume it remains invariant under rigid
translations and rotations. In the case of a Cosserat continuum
they assume invariance of energy balance under
$\mathbf{v}\rightarrow\mathbf{v}+\mathbf{c}_1$ and
$\mathbf{w}\rightarrow\mathbf{w}+\mathbf{c}_2$, where $\mathbf{w}$
is the spatial microstructure velocity field and $\mathbf{c}_1$
and $\mathbf{c}_2$ are arbitrary constant vectors. However, it is
not clear what it means to replace $\mathbf{w}$ by
$\mathbf{w}+\mathbf{c}_2$ in terms of transformations of the
ambient space and microstructure manifolds. In other words, what
group of transformations lead to this replacement and why they
should not affect the macro-velocity field. This seems to be more
or less an assumption convenient for obtaining the desired balance
laws. This assumption leads to conservation of mass and balance of
macro and micro-linear momenta. Then, again they postulate
invariance of local balance of energy under rigid translations and
rotations that transform micro and macro forces simultaneously.
This gives a local form for balance of angular momentum.

\paragraph{The Green-Naghdi-Rivlin Theorem for Structured Continua in
Euclidean Space.}

Let us now study the consequences of postulating invariance of
energy balance under time-dependent isomorphisms of the ambient
Euclidean space with constant velocity for a structured continuum.
Consider balance of energy for
$\varphi_t(\mathcal{U})\subset\varphi_t(\mathcal{B})$ that reads
\begin{equation}
   \frac{d}{dt}\int_{\varphi_t(\mathcal{U})}\rho\left(e+\frac{1}{2}\mathbf{v}\cdot\mathbf{v}\right)dv
   =\int_{\varphi_t(\mathcal{U})} \rho\left(\mathbf{b}\cdot\mathbf{v}+\widetilde{\mathbf{b}}\cdot\widetilde{\mathbf{v}}+r\right)dv
   + \int_{\partial\varphi_t(\mathcal{U})}
   \left(\mathbf{t}\cdot\mathbf{v}+\widetilde{\mathbf{t}}\cdot\widetilde{\mathbf{v}}+h\right)da,
\end{equation}
where for the sake of simplicity, we have ignored the
microstructure inertia. Here $e$ is the internal energy density,
$\mathbf{b}$ is the body force per unit of mass in the deformed
configuration, $\widetilde{\mathbf{b}}$ is the micro-body force
per unit of mass in the deformed configuration, $r$ is heat supply
per unit mass of the deformed configuration, $\mathbf{t}$ is
traction, $\widetilde{\mathbf{t}}$ is micro-traction, and $h$ is
the heat flux. Let us assume that the ambient space is Euclidean,
i.e., $\mathcal{S}=\mathbb{R}^3$. Consider a rigid translation of
the ambient space of the form
\begin{equation}
   \mathbf{x}'=\xi_t(\mathbf{x})=\mathbf{x}+(t-t_0)\mathbf{c},
\end{equation}
where $\mathbf{c}$ is a constant vector field on
$\mathcal{S}=\mathbb{R}^3$. Let us also assume that the director
field
is a vector field on $\mathbb{R}^3$. 
We know that for any $\mathbf{x}\in\mathbb{R}^3$,
$T_{\mathbf{x}}\mathbb{R}^3$ can be identified with $\mathbb{R}^3$
itself. So, we assume that for
$\mathbf{x}=\varphi_t(\mathbf{X})\in\mathbb{R}^3$,
$\mathcal{M}_{\varphi_t(\mathbf{\mathbf{X})}}=T_{\mathbf{x}}\mathbb{R}^3\simeq\mathbb{R}^{3}$.
Note that for a rigid translation of the ambient space
\begin{equation}
   T\xi_t=id,
\end{equation}
where $id$ is the identity map. Therefore, a rigid translation
does not affect the microstructure quantities. Assuming invariance
of balance of energy under arbitrary rigid translations implies
the existence of Cauchy stress and the usual conservation of mass
and balance of energy, i.e.
\begin{eqnarray}
  \dot{\rho}+\rho \operatorname{div}\mathbf{v} &=& 0, \\
  \operatorname{div}\boldsymbol{\sigma}+\rho\mathbf{b} &=&
  \rho\mathbf{a}.
\end{eqnarray}
Next, let us consider a rigid rotation of
$\mathcal{S}=\mathbb{R}^3$ of the form
\begin{equation}
   \mathbf{x}'=\xi_t(\mathbf{x})=e^{\mathbf{\Omega}(t-t_0)}\mathbf{x},
\end{equation}
where $\mathbf{\Omega}$ is a skew-symmetric matrix. Note that
\begin{equation}
   T\xi_t=e^{\mathbf{\Omega}(t-t_0)},~~~TT\xi_t=0.
\end{equation}
We know that
\begin{equation}
   \mathbf{p}'=\xi_{t*}\mathbf{p}=T\xi_t\cdot\mathbf{p}.
\end{equation}
Thus
\begin{equation}
   \widetilde{\mathbf{V}}'=\frac{\partial}{\partial t}\Big|_{\mathbf{X}}\mathbf{p}'=\mathbf{\Omega}e^{\mathbf{\Omega}(t-t_0)}\mathbf{p}
   +e^{\mathbf{\Omega}(t-t_0)}\frac{\partial}{\partial t}\Big|_{\mathbf{X}}\mathbf{p}.
\end{equation}
This means that at $t=t_0$
\begin{equation}
   \widetilde{\mathbf{V}}'=\widetilde{\mathbf{V}}+\mathbf{\Omega}\mathbf{p}.
\end{equation}
Subtracting balance of energy for $\varphi_t(\mathcal{U})$ from
that of $\varphi'_t(\mathcal{U})$ at $t=t_0$, we obtain
\begin{equation}\label{Energy-Rotation}
   \int_{\varphi_t(\mathcal{U})} \rho\mathbf{a}\cdot\mathbf{\Omega}\mathbf{x} ~dv=\int_{\varphi_t(\mathcal{U})}\rho\mathbf{b}\cdot\mathbf{\Omega}\mathbf{x} ~dv
   + \int_{\partial\varphi_t(\mathcal{U})}\mathbf{t}\cdot\mathbf{\Omega}\mathbf{x} ~da
   +\int_{\varphi_t(\mathcal{U})}\rho\widetilde{\mathbf{b}}\cdot\mathbf{\Omega}\mathbf{p} ~dv
   + \int_{\partial\varphi_t(\mathcal{U})}\widetilde{\mathbf{t}}\cdot\mathbf{\Omega}\mathbf{p} ~da.
\end{equation}
We know that
\begin{eqnarray}
  && \label{Surface1} \int_{\partial\varphi_t(\mathcal{U})}\mathbf{t}\cdot\mathbf{\Omega}\mathbf{x}
   ~da=\int_{\varphi_t(\mathcal{U})}\left(\operatorname{div}\boldsymbol{\sigma}\cdot\mathbf{\Omega}\mathbf{x}+
   \boldsymbol{\sigma}:\mathbf{\Omega} \right) dv, \\
  && \label{Surface2} \int_{\partial\varphi_t(\mathcal{U})}\widetilde{\mathbf{t}}\cdot\mathbf{\Omega}\mathbf{p}
   ~da=\int_{\varphi_t(\mathcal{U})}\left[\operatorname{div}\boldsymbol{\widetilde{\sigma}}\otimes\mathbf{p}+
   \boldsymbol{\widetilde{\sigma}}\cdot\nabla\mathbf{p}\right]:\mathbf{\Omega}~dv.
\end{eqnarray}
Substituting (\ref{Surface1}) and (\ref{Surface2}) into
(\ref{Energy-Rotation}) and using the local form of balance of
linear momentum, we obtain
\begin{equation}
   \int_{\varphi_t(\mathcal{U})}\left[\boldsymbol{\sigma}+\operatorname{div}\boldsymbol{\widetilde{\sigma}}\otimes\mathbf{p}+
   \boldsymbol{\widetilde{\sigma}}\cdot\nabla\mathbf{p} \right]:\mathbf{\Omega}~dv=0.
\end{equation}
Because $\mathcal{U}$ is arbitrary, we conclude that
\begin{equation}
   \left[\boldsymbol{\sigma}+\operatorname{div}(\boldsymbol{\widetilde{\sigma}}\otimes\mathbf{p})\right]^{\textsf{T}}
   =\boldsymbol{\sigma}+\operatorname{div}(\boldsymbol{\widetilde{\sigma}}\otimes\mathbf{p}).
\end{equation}
In components this reads as follows:
\begin{equation}
   \sigma^{ab}+\widetilde{\sigma}^{ac}{}_{,c}~p^b+\widetilde{\sigma}^{ac}p^b{}_{,c}=\kappa^{ab}=\kappa^{ba}.
\end{equation}
It is seen that the rigid structure of $\mathbb{R}^3$ and its
isometries does not allow one to obtain a separate balance of
microstructure linear momentum. We will show in the sequel that
when the ambient space is $\mathbb{R}^3$ or, more generally a Riemannian manifold, a generalized covariance can give us such a separate balance of microstructure
linear momentum. We will also see that for a structured continuum
with a scalar microstructure field, e.g., an elastic solid with
distributed voids, one can covariantly obtain a separate scalar
balance of micro-linear momentum.

\vskip 0.3 in
\section {\textbf{A Covariant Theory of Elasticity for Structured Continua with Free Microstructure Manifold}}

In this section we develop a covariant theory of elasticity for
those structured continua for which one can change the spatial and
microstructure frames independently. An example of such continua
is a continuum with voids or a continuum with distributed
``damage", which will be studied in detail in \S5. Let us first
review some important concepts from geometric continuum mechanics.

The reference configuration $\mathcal{B}$ is a submanifold of the
reference configuration manifold $(\mathfrak{B},\mathbf{G})$,
which is a Riemannian manifold. Motion is thought of as an
embedding $\varphi_t:\mathcal{B}\rightarrow \mathcal{S}$, where
$(\mathcal{S},\mathbf{g})$ is the ambient space manifold. An
element $d\mathbf{X}\in T_{\mathbf{X}}\mathcal{B}$ is mapped to
$d\mathbf{x}\in T_{\mathbf{x}}\mathcal{S}$ by the deformation
gradient
\begin{equation}
   d\mathbf{x}=\mathbf{F}\cdot d\mathbf{X}.
\end{equation}
The length of $d\mathbf{x}$ is geometrically important as it
represents the effect of deformation. Note that
\begin{equation}
   \left\langle \! \left\langle d\mathbf{x},d\mathbf{x} \right\rangle \! \right\rangle_{\mathbf{g}}
   =\left\langle \! \left\langle  d\mathbf{X},d\mathbf{X} \right\rangle \! \right\rangle_{\varphi_t^*\mathbf{g}}.
\end{equation}
In this sense the pulled-back metric
$\mathbf{C}=\varphi_t^*\mathbf{g}$ is a measure of deformation.
The material free energy density has the following form
\begin{equation}
   \Psi=\Psi\left(\mathbf{X},\mathbf{F},\mathbf{G},\mathbf{g}\circ \varphi_t\right).
\end{equation}
Let us define the spatial free energy density as
\begin{equation}
   \psi(t,\mathbf{x},\mathbf{g})=\Psi\left(\varphi_t^{-1},\mathbf{F}\circ\varphi_t^{-1},\mathbf{G}\circ\varphi_t^{-1},\mathbf{g}\right).
\end{equation}
Similarly, internal energy density has the following form
\begin{equation}
   e=e(t,\mathbf{x},\mathbf{g}).
\end{equation}
This means that fixing a deformation mapping $\varphi_t$, internal
energy density explicitly depends on time, current position of the
material point and the metric tensor at the current position of
the material point. Note also that $e$ is supported on
$\varphi_t(\mathcal{B})$, i.e. $e=0$ in $\mathcal{S}\setminus
\varphi_t(\mathcal{B})$.

Now let us look at internal energy density for an elastic body
with substructure in which free energy density has the following
form
\begin{equation}
   \Psi=\Psi\left(\mathbf{X},\mathbf{F},\widetilde{\varphi}_t,\widetilde{\mathbf{F}},\mathbf{G},\mathbf{g}\circ \varphi_t
   ,\mathbf{g}_{\mathcal{M}}\circ \widetilde{\varphi}_t\right).
\end{equation}
For a given deformation mapping
$(\varphi_t,\widetilde{\varphi}_t)$ define
\begin{equation}
   \psi(t,\mathbf{x},\mathbf{g},\mathbf{p},\widetilde{\mathbf{g}}_{\mathcal{M}})=\Psi\left(\varphi_t^{-1},
   \mathbf{F}\circ\varphi_t^{-1},\widetilde{\varphi}_t\circ\varphi_t^{-1},\widetilde{\mathbf{F}}\circ\varphi_t^{-1},\mathbf{G}\circ\varphi_t^{-1}
   ,\mathbf{g},\mathbf{p}\circ\varphi_t^{-1},\mathbf{g}_{\mathcal{M}}\circ
   \widetilde{\varphi}_t\circ\varphi_t^{-1}\right),
\end{equation}
where
$\widetilde{\mathbf{g}}_{\mathcal{M}}=\mathbf{g}_{\mathcal{M}}\circ\widetilde{\varphi}\circ\varphi_t^{-1}$.
Similarly, internal energy density has the following form
\begin{equation}
   e=e(t,\mathbf{x},\mathbf{g},\mathbf{p},\widetilde{\mathbf{g}}_{\mathcal{M}}).
\end{equation}
Balance of energy for $\varphi_t(\mathcal{U})\subset \mathcal{S}$
is written as
\begin{eqnarray}
  && \frac{d}{dt}\int_{\varphi_t(\mathcal{U})}\rho(\mathbf{x},t)\left[e(t,\mathbf{x},\mathbf{g},\mathbf{p},\widetilde{\mathbf{g}}_{\mathcal{M}})
   +\frac{1}{2}\left\langle \! \left\langle \mathbf{v},\mathbf{v} \right\rangle \! \right\rangle_{\mathbf{g}} +\kappa
   (\mathbf{p},\widetilde{\mathbf{v}})\right]
    \nonumber \\
  && \label{energy-balance} ~~=\int_{\varphi_t(\mathcal{U})}\rho(\mathbf{x},t)\left(\left\langle \! \left\langle \mathbf{b},\mathbf{v} \right\rangle \! \right\rangle_{\mathbf{g}}
   +\left\langle \!\! \left\langle \widetilde{\mathbf{b}},\widetilde{\mathbf{v}}\right\rangle \!\! \right\rangle_{\widetilde{\mathbf{g}}_{\mathcal{M}}}
   +r\right)+\int_{\partial\varphi_t(\mathcal{U})}\left(\left\langle \! \left\langle \mathbf{t},\mathbf{v} \right\rangle \!
  \right\rangle_{\mathbf{g}}+ \left\langle \!\! \left\langle \widetilde{\mathbf{t}},\widetilde{\mathbf{v}} \right\rangle
  \!\!
  \right\rangle_{\widetilde{\mathbf{g}}_{\mathcal{M}}}+h\right)da,
\end{eqnarray}
where we think of $\rho(\mathbf{x},t)$ as a $3$-form and
$\widetilde{\mathbf{b}}$ and $\widetilde{\mathbf{t}}$ are
microstructure body force and traction vector fields,
respectively. For the sake of simplicity, let us assume that the
microstructure kinetic energy has the following form
\begin{equation}
   \kappa (\mathbf{p},\widetilde{\mathbf{v}})=\frac{1}{2}j\left\langle \! \left\langle \widetilde{\mathbf{v}},\widetilde{\mathbf{v}} \right\rangle \!
   \right\rangle_{\widetilde{\mathbf{g}}_{\mathcal{M}}},
\end{equation}
where we assume the microstructure inertia $j$ is a scalar.

All the physical processes happen in $\mathcal{S}$ and thus
balance of energy is written on subsets of
$\varphi_t(\mathcal{B})\subset\mathcal{S}$. Standard traction is a
vector field on $\mathcal{S}$ and the microstructure traction is a
vector field on $\mathcal{M}$. The standard and microstructure
tractions have the following coordinate representations
\begin{equation}
   \mathbf{t}(\mathbf{x},t)=t^a\mathbf{e}_a,~~~~~\widetilde{\mathbf{t}}(\mathbf{x},t)=\widetilde{t}^{\alpha}~\widetilde{\mathbf{e}}_{\alpha},
\end{equation}
where $\{\mathbf{e}_a\}$ and $\{\widetilde{\mathbf{e}}_a\}$ are
bases for $T_{\mathbf{x}}\mathcal{S}$ and
$T_{\mathbf{p}}\mathcal{M}$, respectively. Similarly, the stress
tensors have the following local representations
\begin{equation}
   \boldsymbol{\sigma}(\mathbf{x},t)=\sigma^{ab}~\mathbf{e}_a\otimes\mathbf{e}_b,~~~~~\widetilde{\boldsymbol{\sigma}}(\mathbf{x},t)
   =\widetilde{\sigma}^{\alpha b}~\widetilde{\mathbf{e}}_{\alpha}\otimes\mathbf{e}_{b}.
\end{equation}
The first Piola Kirchhoff stresses for the standard deformation
and the microstructure deformation are obtained by the following
Piola transformations
\begin{equation}
   P^{aA}=J(\mathbf{F}^{-1})^A{}_{b}~\sigma^{ab},~~~~~\widetilde{P}^{\alpha A}=J(\mathbf{F}^{-1})^A{}_{b}~\widetilde{\sigma}^{\alpha b},
\end{equation}
where
$J=\sqrt{\frac{\textrm{det}\mathbf{g}}{\textrm{det}\mathbf{G}}}~\textrm{det}\mathbf{F}$.
These transformations ensure that
\begin{equation}
   \mathbf{t}~da=\mathbf{T}~dA~~~~~\textrm{and}~~~~~\widetilde{\mathbf{t}}~da=\widetilde{\mathbf{T}}~dA.
\end{equation}
Now this means that in terms of contributions of tractions to
balance of energy we have
\begin{equation}
   \left\langle \! \left\langle\mathbf{t},\mathbf{v}\right\rangle \! \right\rangle_{\mathbf{g}} da
   =\left\langle \! \left\langle \mathbf{T},\mathbf{V}\right\rangle \! \right\rangle_{\mathbf{g}}dA
   ~~~~~\textrm{and}~~~~~
   \left\langle \!\! \left\langle\widetilde{\mathbf{t}},\widetilde{\mathbf{v}}\right\rangle \!\! \right\rangle_{\mathbf{g}_{\mathcal{M}}} da
   =\left\langle \!\! \left\langle \widetilde{\mathbf{T}},\widetilde{\mathbf{V}}\right\rangle \!\! \right\rangle_{\mathbf{g}_{\mathcal{M}}}dA.
\end{equation}
For $\mathcal{U}\subset\mathcal{B}$, material energy balance can
be written as
\begin{eqnarray}\label{material-energy-balance}
  && \frac{d}{dt}\int_{\mathcal{U}}\rho_0(\mathbf{X},t)\left[E(t,\mathbf{X},\mathbf{g},\mathbf{g}_{\mathcal{M}})
   +\frac{1}{2}\left\langle \! \left\langle \mathbf{V},\mathbf{V} \right\rangle \! \right\rangle_{\mathbf{g}}
   +\frac{1}{2}J\left\langle \!\! \left\langle \widetilde{\mathbf{V}},\widetilde{\mathbf{V}}\right\rangle \!\! \right\rangle_{\mathbf{g}_{\mathcal{M}}}
   \right]
    \nonumber \\
  && \label{material-energy-balance}~~=\int_{\mathcal{U}}\rho_0(\mathbf{X},t)\left(\left\langle \! \left\langle \mathbf{B},\mathbf{V} \right\rangle \! \right\rangle_{\mathbf{g}}
   +\left\langle \!\! \left\langle \widetilde{\mathbf{B}},\widetilde{\mathbf{V}}\right\rangle \!\! \right\rangle_{\widetilde{\mathbf{g}}_{\mathcal{M}}}
   +R\right)+\int_{\partial\mathcal{U}}\left(\left\langle \! \left\langle \mathbf{T},\mathbf{V} \right\rangle \!
  \right\rangle_{\mathbf{g}}
  + \left\langle \!\! \left\langle \widetilde{\mathbf{T}},\widetilde{\mathbf{V}}\right\rangle \!\!
  \right\rangle_{\mathbf{g}_{\mathcal{M}}}+H\right)dA,
\end{eqnarray}
where again $\rho_0$ is a $3$-form.

\subsection{\textbf{Covariance of Energy Balance}}
\label{sec:cov}

Let us assume that for each $\mathbf{x}\in\mathcal{S}$, the
microstructure manifold is completely independent of
$\mathcal{S}$. In other words, a change of frame in
$\mathcal{S}(\textrm{or}~\mathcal{M})$ does not affect
$\mathcal{M}(\textrm{or}~\mathcal{S})$ and quantities defined on
it. An example of a structured continuum with this type of
microstructure manifold is a structured continuum with a scalar
director field, although there are other possibilities. We show in
this subsection that postulating energy balance and its invariance
under time-dependent changes of frame in $\mathcal{S}$ and
$\mathcal{M}$ results in conservation of mass and micro-inertia,
two balances of linear and angular momenta, and two Doyle-Ericksen
formulas, one for the Cauchy stress and one for the micro-Cauchy
stress.

\vskip 0.1 in {\theorem If balance of energy holds and if it is
invariant under arbitrary spatial and microstructure
diffeomorhisms $\xi_t:\mathcal{S}\rightarrow\mathcal{S}$ and
$\eta_t:\mathcal{M}\rightarrow\mathcal{M}$, then there exist
second-order tensors $\boldsymbol{\sigma}$ and
$\widetilde{\boldsymbol{\sigma}}$ such that
\begin{equation}
   \mathbf{t}=\left\langle \! \left\langle \boldsymbol{\sigma},\mathbf{n}   \right\rangle \! \right\rangle_{\mathbf{g}}
   ~~~\textrm{and}~~~
   \widetilde{\mathbf{t}}=\left\langle \! \left\langle \widetilde{\boldsymbol{\sigma}},\mathbf{n}   \right\rangle \!
   \right\rangle_{\mathbf{g}},
\end{equation}
and
\begin{eqnarray}
  && \mathbf{L}_{\mathbf{v}}\rho=0, \\
  && \mathbf{L}_{\mathbf{v}}j=0, \\
  && \operatorname{div}\boldsymbol{\sigma}+\rho\mathbf{b}=\rho \mathbf{a}, \\
  && \operatorname{div}\widetilde{\boldsymbol{\sigma}}+\rho\widetilde{\mathbf{b}}=\rho j\widetilde{\mathbf{a}}, \\
  && \boldsymbol{\sigma}=\boldsymbol{\sigma}^{\textsf{T}}, \\
  && \left(\mathbf{F}_0\widetilde{\boldsymbol{\sigma}}\right)^{\textsf{T}}=\mathbf{F}_0\widetilde{\boldsymbol{\sigma}}, \\
  && 2\rho\frac{\partial e}{\partial\mathbf{g}}=\boldsymbol{\sigma}, \\
  && \mathbf{F}_0\widetilde{\boldsymbol{\sigma}}=2\rho\frac{\partial e}{\partial
  \widetilde{\mathbf{g}}_{\mathcal{M}}},
\end{eqnarray}
where $\operatorname{div}$ is divergence with respect to the
metric $\mathbf{g}$,
$\mathbf{F}_0=\widetilde{\mathbf{F}}\mathbf{F}^{-1}$ and $\eta_t$
acts on all the microstructure fibers simultaneously.}

\begin{figure}[hbt]
\vspace*{0.3in}
\begin{center}
\includegraphics[scale=0.7,angle=0]{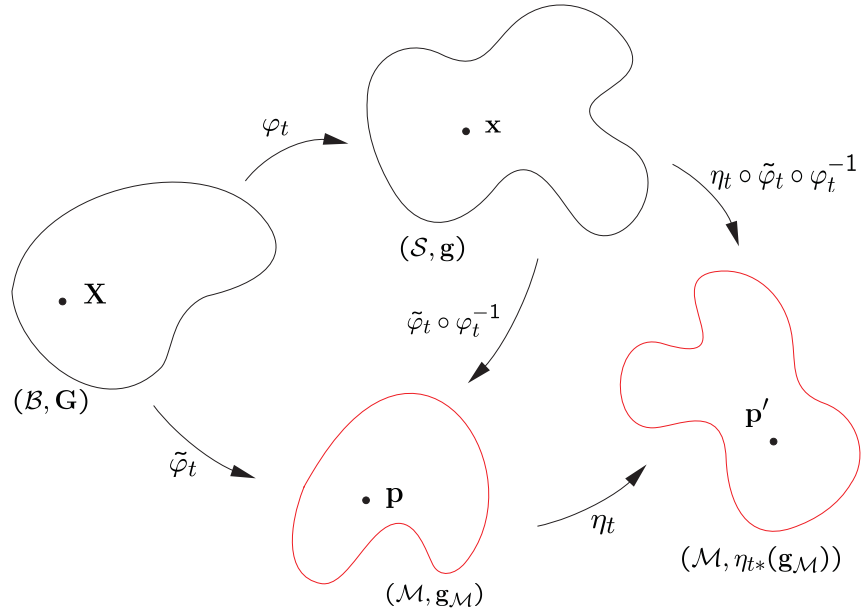}
\end{center}
\caption{\footnotesize A microstructure change of frame.}
\label{Cosserat_Framing1}
\end{figure}

\paragraph{Proof:} Let us consider spatial and microstructure diffeomorphisms separately.

\paragraph{Microstructure covariance of energy balance.} Consider a microstructure diffeomorphism
$\eta_t:\mathcal{M}\rightarrow\mathcal{M}$ (see Fig.
\ref{Cosserat_Framing1}) and assume that
\begin{equation}
   \eta_t\big|_{t=t_0}=id.
\end{equation}
Invariance of energy balance under
$\eta_t:\mathcal{M}\rightarrow\mathcal{M}$ means that balance of
energy in the new frame has the following form
\begin{eqnarray}
  && \frac{d}{dt}\int_{\varphi_t(\mathcal{U})}\rho(\mathbf{x},t)\left[e'(t,\mathbf{x},\mathbf{g},\mathbf{p}',\widetilde{\mathbf{g}}_{\mathcal{M}})
   +\frac{1}{2}\left\langle \! \left\langle \mathbf{v},\mathbf{v} \right\rangle \! \right\rangle_{\mathbf{g}}
   +\frac{1}{2}j'\left\langle \! \left\langle \widetilde{\mathbf{v}}',\widetilde{\mathbf{v}}'\right\rangle \! \right\rangle_{\widetilde{\mathbf{g}}_{\mathcal{M}}}
   \right]
    \nonumber \\
  && ~~=\int_{\varphi_t(\mathcal{U})}\rho(\mathbf{x},t)\left(\left\langle \! \left\langle \mathbf{b},\mathbf{v} \right\rangle \! \right\rangle_{\mathbf{g}}
   +\left\langle \!\! \left\langle \widetilde{\mathbf{b}}',\widetilde{\mathbf{v}}'\right\rangle \!\! \right\rangle_{\widetilde{\mathbf{g}}_{\mathcal{M}}}
   +r\right)+\int_{\partial\varphi_t(\mathcal{U})}\left(\left\langle \! \left\langle \mathbf{t},\mathbf{v} \right\rangle \!
  \right\rangle_{\mathbf{g}}+ \left\langle \!\! \left\langle \widetilde{\mathbf{t}}',\widetilde{\mathbf{v}}' \right\rangle
  \!\!
  \right\rangle_{\widetilde{\mathbf{g}}_{\mathcal{M}}}+h\right)da.
\end{eqnarray}
Note that
\begin{equation}
   e'(t,\mathbf{x},\mathbf{g},\mathbf{p}',\widetilde{\mathbf{g}}_{\mathcal{M}})=e(t,\mathbf{x},\mathbf{g},\mathbf{p},\eta_t^{*}\widetilde{\mathbf{g}}_{\mathcal{M}}).
\end{equation}
Thus
\begin{equation}
   \frac{d}{dt}\Big|_{t=t_0}=\dot{e}+\frac{\partial e}{\partial
   \widetilde{\mathbf{g}}_{\mathcal{M}}}:\mathfrak{L}_{\mathbf{z}}\widetilde{\mathbf{g}}_{\mathcal{M}},
\end{equation}
where
\begin{equation}
   \mathbf{z}=\frac{\partial}{\partial t}\Big|_{t=t_0}\eta_t.
\end{equation}
Note also that
\begin{equation}
   \widetilde{\mathbf{v}}'\big|_{t=t_0}=\widetilde{\mathbf{v}}+\mathbf{z}.
\end{equation}
Assuming that
$\widetilde{\mathbf{b}}'-j'\widetilde{\mathbf{a}}'=\eta_{t*}(\widetilde{\mathbf{b}}-j\widetilde{\mathbf{a}})$,
at $t=t_0$ we obtain
\begin{eqnarray}
  && \int_{\varphi_t(\mathcal{U})}\mathbf{L}_{\mathbf{v}}\rho\left(e+\left\langle \! \left\langle \mathbf{v},\mathbf{v} \right\rangle \! \right\rangle_{\mathbf{g}}
  +\frac{1}{2}j\left\langle \! \left\langle \widetilde{\mathbf{v}}+\mathbf{z},\widetilde{\mathbf{v}}+\mathbf{z}\right\rangle \! \right\rangle_{\widetilde{\mathbf{g}}_{\mathcal{M}}}\right)
     \nonumber \\
  && ~~+\int_{\varphi_t(\mathcal{U})}\rho\left(\dot{e}+\frac{\partial e}{\partial
   \widetilde{\mathbf{g}}_{\mathcal{M}}}:\mathfrak{L}_{\mathbf{z}}\widetilde{\mathbf{g}}_{\mathcal{M}}
   +j\left\langle \! \left\langle \widetilde{\mathbf{a}},\mathbf{z}\right\rangle \! \right\rangle_{\widetilde{\mathbf{g}}_{\mathcal{M}}}
   +\frac{1}{2}\mathbf{L}_{\mathbf{v}}j\left\langle \! \left\langle \widetilde{\mathbf{v}}+\mathbf{z},\widetilde{\mathbf{v}}+\mathbf{z}\right\rangle \! \right\rangle_{\widetilde{\mathbf{g}}_{\mathcal{M}}}  \right)    \nonumber \\
  && ~~=\int_{\varphi_t(\mathcal{U})}\rho\left(\left\langle \! \left\langle \mathbf{b},\mathbf{v} \right\rangle \! \right\rangle_{\mathbf{g}}
   +\left\langle \!\! \left\langle \widetilde{\mathbf{b}},\widetilde{\mathbf{v}}+\mathbf{z}\right\rangle \!\! \right\rangle_{\widetilde{\mathbf{g}}_{\mathcal{M}}}
   +r\right)+\int_{\partial\varphi_t(\mathcal{U})}\left(\left\langle \! \left\langle \mathbf{t},\mathbf{v} \right\rangle \!
  \right\rangle_{\mathbf{g}}+ \left\langle \!\! \left\langle \widetilde{\mathbf{t}},\widetilde{\mathbf{v}}+\mathbf{z} \right\rangle \!\!
  \right\rangle_{\widetilde{\mathbf{g}}_{\mathcal{M}}}+h\right)da.
\end{eqnarray}
Replacing $\rho$ by $\rho dv$ and subtracting balance of energy
(\ref{energy-balance}) from the above identity and considering the
fact that $\mathbf{z}$ and $\mathcal{U}$ are arbitrary, one
obtains
\begin{eqnarray}
  && \label{mass-mass-conservation} \mathbf{L}_{\mathbf{v}}(\rho j)=0, \\
  && \nonumber \\
  && \label{Cauchy's Theorem} \int_{\varphi_t(\mathcal{U})}\rho\frac{\partial e}{\partial
   \widetilde{\mathbf{g}}_{\mathcal{M}}}:\mathfrak{L}_{\mathbf{z}}\widetilde{\mathbf{g}}_{\mathcal{M}}~dv
   =\int_{\varphi_t(\mathcal{U})}\rho\left\langle \!\! \left\langle \widetilde{\mathbf{b}},\mathbf{z}\right\rangle \!\!
   \right\rangle_{\widetilde{\mathbf{g}}_{\mathcal{M}}}dv
   +\int_{\partial\varphi_t(\mathcal{U})}\left\langle \!\! \left\langle \widetilde{\mathbf{t}},\mathbf{z} \right\rangle \!\! \right\rangle_{\widetilde{\mathbf{g}}_{\mathcal{M}}}da.
\end{eqnarray}
Applying Cauchy's theorem (see \cite{MaHu1983}) to \eqref{Cauchy's
Theorem}, one concludes that there exists a second-order tensor
$\widetilde{\boldsymbol{\sigma}}$ such that
\begin{equation}\label{OrderParameterCT}
   \widetilde{\mathbf{t}}=\left\langle \! \left\langle \widetilde{\boldsymbol{\sigma}},\mathbf{n}   \right\rangle \!
   \right\rangle_{\mathbf{g}}.
\end{equation}
Now let us simplify the surface integral. \vskip 0.1 in {\lemma
The contribution of microstructure traction has the following
simplified form.
\begin{equation}
   \int_{\partial\varphi_t(\mathcal{U})}\left\langle \!\! \left\langle \widetilde{\mathbf{t}},\mathbf{z} \right\rangle \!\!
   \right\rangle_{\widetilde{\mathbf{g}}_{\mathcal{M}}}da
   =\int_{\varphi_t(\mathcal{U})}\left[\left\langle \! \left\langle  \operatorname{div}\widetilde{\boldsymbol{\sigma}}  ,\mathbf{z}  \right\rangle \! \right\rangle_{\widetilde{\mathbf{g}}_{\mathcal{M}}}
   +\mathbf{F}_0\widetilde{\boldsymbol{\sigma}}:\frac{1}{2}\mathfrak{L}_{\mathbf{z}}\widetilde{\mathbf{g}}_{\mathcal{M}}
   +\mathbf{F}_0\widetilde{\boldsymbol{\sigma}}:\boldsymbol{\omega}_{\mathcal{M}} \right]dv.
\end{equation}
}

\paragraph{Proof:}
\begin{equation}
   \int_{\partial\varphi_t(\mathcal{U})}\left\langle \!\! \left\langle \widetilde{\mathbf{t}},\mathbf{z} \right\rangle \!\! \right\rangle_{\widetilde{\mathbf{g}}_{\mathcal{M}}}
   =\int_{\partial\varphi_t(\mathcal{U})}\sigma^{\alpha b}n^cg_{bc}z^{\beta}(g_{\mathcal{M}})_{\alpha\beta}~da
   =\int_{\varphi_t(\mathcal{U})}\left[\sigma^{\alpha b}z^{\beta}(g_{\mathcal{M}})_{\alpha\beta}\right]_{|b}dv.
\end{equation}
But because
$(g_{\mathcal{M}})_{{\alpha\beta|b}}=(g_{\mathcal{M}})_{{\alpha\beta|\gamma}}(F_0)^{\gamma}{}_{b}=0$,
we have
\begin{equation}
    \left[\sigma^{\alpha b}z^{\beta}(g_{\mathcal{M}})_{\alpha\beta}\right]_{|b} = \left[\sigma^{\alpha b}z^{\beta}\right]_{|b}(g_{\mathcal{M}})_{\alpha\beta}
   = \sigma^{\alpha b}{}_{|b}z^{\beta}(g_{\mathcal{M}})_{\alpha\beta}+z^{\beta}{}_{|b}\sigma^{\alpha
   b}(g_{\mathcal{M}})_{\alpha\beta}.
\end{equation}
Note that
\begin{equation}
   z^{\beta}{}_{|b}(g_{\mathcal{M}})_{\alpha\beta}=z_{\alpha|\gamma}\left(F_0\right)^{\lambda}{}_{b}.~~~~~\square
\end{equation}
Now, because $\mathbf{z}$ and $\mathcal{U}$ are arbitrary from
(\ref{Cauchy's Theorem}) one obtains
\begin{eqnarray}
  && \label{DE-M}\mathbf{F}_0\widetilde{\boldsymbol{\sigma}}=2\rho\frac{\partial e}{\partial \widetilde{\mathbf{g}}_{\mathcal{M}}}, \\
  && \left(\mathbf{F}_0\widetilde{\boldsymbol{\sigma}}\right)^{\textsf{T}}=\mathbf{F}_0\widetilde{\boldsymbol{\sigma}}, \\
  &&
  \operatorname{div}\widetilde{\boldsymbol{\sigma}}+\rho\widetilde{\mathbf{b}}=\rho j\widetilde{\mathbf{a}}.
\end{eqnarray}

\begin{figure}[hbt]
\vspace*{0.3in}
\begin{center}
\includegraphics[scale=0.7,angle=0]{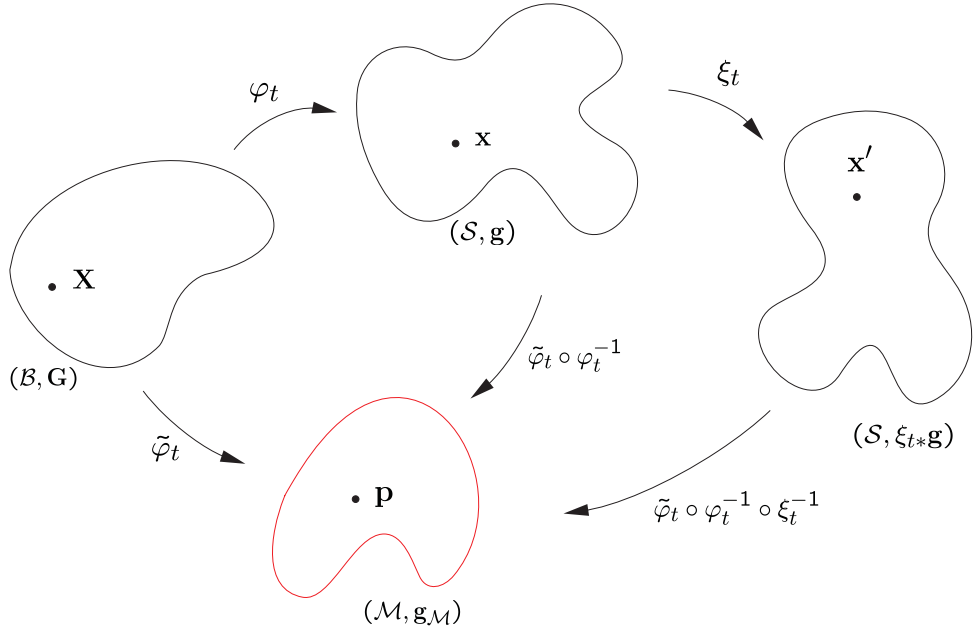}
\end{center}
\caption{\footnotesize A spatial change of frame in a continuum
with microstructure.} \label{material-framing}
\end{figure}

\paragraph{Spatial covariance of energy balance.} Invariance of energy balance under an arbitrary diffeomorphism
$\xi_t:\mathcal{S}\rightarrow\mathcal{S}$ means that (see Fig.
\ref{material-framing})
\begin{eqnarray}
   && \frac{d}{dt}\int_{\varphi'_t(\mathcal{U})}\rho'(\mathbf{x}',t)\left[e'(t,\mathbf{x}',\mathbf{g},\mathbf{g}_{\mathcal{M}})
   +\frac{1}{2}\left\langle \! \left\langle \mathbf{v}',\mathbf{v}' \right\rangle \! \right\rangle_{\mathbf{g}}
      +\frac{1}{2}j'\left\langle \! \left\langle \widetilde{\mathbf{v}}',\widetilde{\mathbf{v}}'\right\rangle \! \right\rangle_{\widetilde{\mathbf{g}}_{\mathcal{M}}}\right] \nonumber \\
   && ~=\int_{\varphi'_t(\mathcal{U})}\rho'(\mathbf{x}',t)\left(\left\langle \! \left\langle \mathbf{b}',\mathbf{v}' \right\rangle \! \right\rangle_{\mathbf{g}}
    +\left\langle \!\! \left\langle \widetilde{\mathbf{b}}',\widetilde{\mathbf{v}}' \right\rangle \!\! \right\rangle_{\widetilde{\mathbf{g}}_{\mathcal{M}}} \!\!\!+r'
    \right)+\int_{\partial\varphi'_t(\mathcal{U})}\left(\left\langle \! \left\langle \mathbf{t}',\mathbf{v}' \right\rangle \!
  \right\rangle_{\mathbf{g}}+ \left\langle \!\! \left\langle \widetilde{\mathbf{t}}',\widetilde{\mathbf{v}}' \right\rangle \!\! \right\rangle_{\widetilde{\mathbf{g}}_{\mathcal{M}}}\!\!\!+h'\right)da',
\end{eqnarray}
where $\varphi'_t=\xi_t\circ\varphi_t$. We also assume that
\begin{equation}
   \xi_t\big|_{t=t_0}=id.
\end{equation}
The relation between primed and unprimed quantities are dictated
by Cartan's spacetime theory, i.e.
\begin{equation}
    \rho'(\mathbf{x}',t)=\xi_{*}\rho(\mathbf{x},t),~\mathbf{t}'=\xi_{*}\mathbf{t},~\widetilde{\mathbf{t}}'=\xi_{*}\widetilde{\mathbf{t}},
    ~r'(\mathbf{x}',t)=r(\mathbf{x},t),~ h'(\mathbf{x}',t)=h(\mathbf{x},t).
\end{equation}
The internal energy density has the following transformation
\begin{equation}
   e'(t,\mathbf{x}',\mathbf{g},\widetilde{\mathbf{g}}_{\mathcal{M}})=e(t,\mathbf{x},\xi^{*}\mathbf{g},\mathbf{p},\widetilde{\mathbf{g}}_{\mathcal{M}}).
\end{equation}
Thus
\begin{equation}
   \frac{d}{dt}\Big|_{t=t_0}e'=\dot{e}+\frac{\partial e}{\partial
   \mathbf{g}}:\mathfrak{L}_{\mathbf{w}}\mathbf{g},
\end{equation}
where
\begin{equation}
   \mathbf{w}=\frac{\partial}{\partial t}\Big|_{t=t_0}\xi_t.
\end{equation}
Spatial velocity has the following transformation
\begin{equation}
   \mathbf{v}'=\xi_{*}\mathbf{v}+\mathbf{w}_t.
\end{equation}
Thus, at $t=t_0$, $\mathbf{v}'=\mathbf{v}+\mathbf{w}$. Also
\begin{equation}
   \widetilde{\mathbf{v}}'=\widetilde{\mathbf{V}}\circ\varphi_t^{-1}\circ\xi_t^{-1}=\widetilde{\mathbf{v}}\circ\xi_t^{-1}.
\end{equation}
Therefore, at $t=t_0$
\begin{equation}
   \widetilde{\mathbf{v}}'=\widetilde{\mathbf{v}}.
\end{equation}
Assuming that
$\mathbf{b}'-\mathbf{a}'=\xi_{t*}(\mathbf{b}-\mathbf{a})$
\citep{MaHu1983} and noting that
$\widetilde{\mathbf{b}}'-\widetilde{\mathbf{a}}'=\widetilde{\mathbf{b}}-\widetilde{\mathbf{a}}$,
balance of energy in the new frame at $t=t_0$ reads
\begin{eqnarray}
  && \int_{\varphi_t(\mathcal{U})}\mathbf{L}_{\mathbf{v}}\rho\left(e+\frac{1}{2}\left\langle \! \left\langle \mathbf{v}+\mathbf{w},\mathbf{v}+\mathbf{w} \right\rangle \! \right\rangle_{\mathbf{g}}
  +\frac{1}{2}j\left\langle \! \left\langle \widetilde{\mathbf{v}},\widetilde{\mathbf{v}} \right\rangle \!
  \right\rangle_{\widetilde{\mathbf{g}}_{\mathcal{M}}}\right)    \nonumber \\
  && ~~~+\int_{\varphi_t(\mathcal{U})}\rho\left(\dot{e}+\frac{\partial e}{\partial \mathbf{g}}:\mathfrak{L}_{\mathbf{w}}\mathbf{g}+\left\langle \! \left\langle \mathbf{v}+\mathbf{w},\mathbf{a} \right\rangle \! \right\rangle_{\mathbf{g}}
  +j\left\langle \! \left\langle \widetilde{\mathbf{v}},\widetilde{\mathbf{a}} \right\rangle \!
  \right\rangle_{\widetilde{\mathbf{g}}_{\mathcal{M}}}
  +\frac{1}{2}\mathbf{L}_{\mathbf{v}}j\left\langle \! \left\langle \widetilde{\mathbf{v}},\widetilde{\mathbf{v}} \right\rangle \! \right\rangle_{\widetilde{\mathbf{g}}_{\mathcal{M}}}\right)
    \nonumber \\
  && \label{energy-balance-framing} ~=\int_{\varphi_t(\mathcal{U})}\rho\left(\left\langle \! \left\langle \mathbf{b},\mathbf{v}+\mathbf{w} \right\rangle \! \right\rangle_{\mathbf{g}}
   +\left\langle \!\! \left\langle \widetilde{\mathbf{b}},\widetilde{\mathbf{v}}\right\rangle \!\! \right\rangle_{\widetilde{\mathbf{g}}_{\mathcal{M}}}
   +r\right)+\int_{\partial\varphi_t(\mathcal{U})}\left(\left\langle \! \left\langle \mathbf{t},\mathbf{v}+\mathbf{w} \right\rangle \!
  \right\rangle_{\mathbf{g}}+ \left\langle \!\! \left\langle \widetilde{\mathbf{t}},\widetilde{\mathbf{v}} \right\rangle \!\!
  \right\rangle_{\widetilde{\mathbf{g}}_{\mathcal{M}}}+h\right)da.
\end{eqnarray}
Subtracting (\ref{energy-balance}) from
(\ref{energy-balance-framing}) and considering the fact that
$\mathbf{w}$ and $\mathcal{U}$ are arbitrary, we obtain
conservation of mass $\mathbf{L}_{\mathbf{v}}\rho=0$ and using it
in (\ref{mass-mass-conservation}) we obtain balance of
microstructure inertia
\begin{equation}
   \mathbf{L}_{\mathbf{v}}j=0.
\end{equation}
Now using conservation of mass and microstructure inertia, and
replacing $\rho$ by $\rho dv$ in (\ref{energy-balance-framing}),
one obtains
\begin{equation}\label{energy-balance-framing1}
    \int_{\varphi_t(\mathcal{U})}\rho\left(\frac{\partial e}{\partial \mathbf{g}}:\mathfrak{L}_{\mathbf{w}}\mathbf{g}+\left\langle\!\left\langle \mathbf{w},\mathbf{a} \right\rangle \! \right\rangle_{\mathbf{g}}
  \right)dv=\int_{\varphi_t(\mathcal{U})}\rho\left(\left\langle\!\left\langle \mathbf{b},\mathbf{w} \right\rangle \! \right\rangle_{\mathbf{g}}\right)dv
+\int_{\partial\varphi_t(\mathcal{U})}\left(\left\langle \!
\left\langle \mathbf{t},\mathbf{w} \right\rangle \!
  \right\rangle_{\mathbf{g}}\right)da.
\end{equation}
Applying Cauchy's theorem to the above identity and considering
(\ref{OrderParameterCT}) shows that there exists a second-order
tensor $\boldsymbol{\sigma}$ such that
\begin{equation}
   \mathbf{t}=\left\langle \! \left\langle \boldsymbol{\sigma},\mathbf{n}   \right\rangle \! \right\rangle_{\mathbf{g}}.
\end{equation}
Now let us look at the surface integral in
(\ref{energy-balance-framing1}). This surface integral is
simplified to read
\begin{equation}\label{surface}
    \int_{\partial\varphi_t(\mathcal{U})}\left\langle \! \left\langle \mathbf{t},\mathbf{w} \right\rangle \!
  \right\rangle_{\mathbf{g}} da = \int_{\varphi_t(\mathcal{U})}\left\langle \!\left\langle \operatorname{div}\boldsymbol{\sigma}
  ,\mathbf{w} \right\rangle \! \right\rangle_{\mathbf{g}}dv+\int_{\varphi_t(\mathcal{U})}\left(\boldsymbol{\sigma}:\frac{1}{2}\mathfrak{L}_{\mathbf{w}}\mathbf{g}
   +\boldsymbol{\sigma}:\boldsymbol{\omega} \right)dv,
\end{equation}
where $\boldsymbol{\omega}$ has the coordinate representation
$\omega_{ab}=\frac{1}{2}(w_{a|b}-w_{b|a})$. Substituting
(\ref{surface}) into (\ref{energy-balance-framing1}) yields
\begin{equation}
    \int_{\varphi_t(\mathcal{U})}\left(2\rho\frac{\partial e}{\partial\mathbf{g}}-\boldsymbol{\sigma}\right):\frac{1}{2}\mathfrak{L}_{\mathbf{w}}\mathbf{g}~dv
  + \int_{\varphi_t(\mathcal{U})} \boldsymbol{\sigma}:\boldsymbol{\omega}~dv-\int_{\varphi_t(\mathcal{U})} \left\langle \!\left\langle \operatorname{div}\boldsymbol{\sigma}+\rho\left(\mathbf{b}
  -\mathbf{a}\right),\mathbf{w} \right\rangle \! \right\rangle_{\mathbf{g}} dv=0.
\end{equation}
Because $\mathcal{U}$ and $\mathbf{w}$ are arbitrary we conclude
that
\begin{eqnarray}
  && \label{DE-Energy} 2\rho\frac{\partial e}{\partial\mathbf{g}}=\boldsymbol{\sigma},  \\
  && \boldsymbol{\sigma}
  =\boldsymbol{\sigma}^{\textsf{T}}, \\
  && \operatorname{div}\boldsymbol{\sigma}+\rho\mathbf{b}=\rho
  \mathbf{a}.~~~~~~\square
\end{eqnarray}
Next, we study the effect of material diffeomorphisms on balance
of energy.

\vskip 0.4 in \subsection{\textbf{Transformation of Energy Balance
under Material Diffeomorphisms}}

It was shown in \cite{YaMaOr2006} that, in general, energy
balance cannot be invariant under diffeomorphisms of the reference
configuration and what one should be looking for instead is the
way in which energy balance transforms under material
diffeomorphisms. In this subsection we first obtain such a
transformation formula for a continuum with microstructure under
an arbitrary time-dependent material diffeomorphism (see Eq.
\eqref{Energy-Balance-Transformation}) and then obtain the
conditions under which balance of energy can be materially
covariant.

\paragraph{\textbf{The Material Energy Balance Transformation Formula.}}
Let us begin with a discussion of how energy balance transforms
under material diffeomorphisms. Let us define
\begin{equation}
    E(t,\mathbf{X},\mathbf{G})=
    E \left(\mathbf{X},\mathbf{F}(\mathbf{X}),\widetilde{\varphi}_t(\mathbf{X}),\widetilde{\mathbf{F}}(\mathbf{X}),\mathbf{g}(\varphi_t(\mathbf{X})),\mathbf{g}_{\mathcal{M}}(\widetilde{\varphi}_t(\mathbf{X})),\mathbf{G}\right),
\end{equation}
where $E$ is the material internal energy density per unit of
undeformed mass. Material (Lagrangian) energy balance
(\ref{material-energy-balance}) can be simplified to read
\begin{eqnarray}
  && \int_{\mathcal{U}}\frac{d}{dt}\left[\rho_0\left(E(t,\mathbf{X},\mathbf{G})
   +\frac{1}{2}\left\langle \! \left\langle \mathbf{V},\mathbf{V} \right\rangle \! \right\rangle_{\mathbf{g}}
   +\frac{1}{2}J\left\langle \!\! \left\langle \widetilde{\mathbf{V}},\widetilde{\mathbf{V}}\right\rangle \!\! \right\rangle_{\mathbf{g}_{\mathcal{M}}}
   \right)\right]
    \nonumber \\
  && ~~=\int_{\mathcal{U}}\rho_0\left(\left\langle \! \left\langle \mathbf{B},\mathbf{V} \right\rangle \! \right\rangle_{\mathbf{g}}
   +\left\langle \!\! \left\langle \widetilde{\mathbf{B}},\widetilde{\mathbf{V}}\right\rangle \!\! \right\rangle_{\widetilde{\mathbf{g}}_{\mathcal{M}}}
   +R\right)+\int_{\partial\mathcal{U}}\left(\left\langle \! \left\langle \mathbf{T},\mathbf{V} \right\rangle \!
  \right\rangle_{\mathbf{g}}
  + \left\langle \!\! \left\langle \widetilde{\mathbf{T}},\widetilde{\mathbf{V}}\right\rangle \!\!
  \right\rangle_{\mathbf{g}_{\mathcal{M}}}+H\right)dA,
\end{eqnarray}
where $\mathcal{U}$ is an arbitrary nice subset of the reference
configuration $\mathcal{B}$, $\mathbf{B}$ and
$\widetilde{\mathbf{B}}$ are body force and microstructure body
force, respectively, per unit undeformed mass,
$\mathbf{V}(\mathbf{X},t)$ and
$\widetilde{\mathbf{V}}(\mathbf{X},t)$ are the material velocity
and microstructure material velocity, respectively,
$\rho_0(\mathbf{X},t)$ is the material density, $R(\mathbf{X},t)$
is the heat supply per unit undeformed mass, and
$H(\mathbf{X},t,\hat{\mathbf{N}})$ is the heat flux across a
surface with normal $\hat{\mathbf{N}}$ in the undeformed
configuration (normal to $\partial\mathcal{U}$ at
$\mathbf{X}\in\partial\mathcal{U}$).

\begin{figure}[hbt]
\vspace*{0.3in}
\begin{center}
\includegraphics[scale=0.75,angle=0]{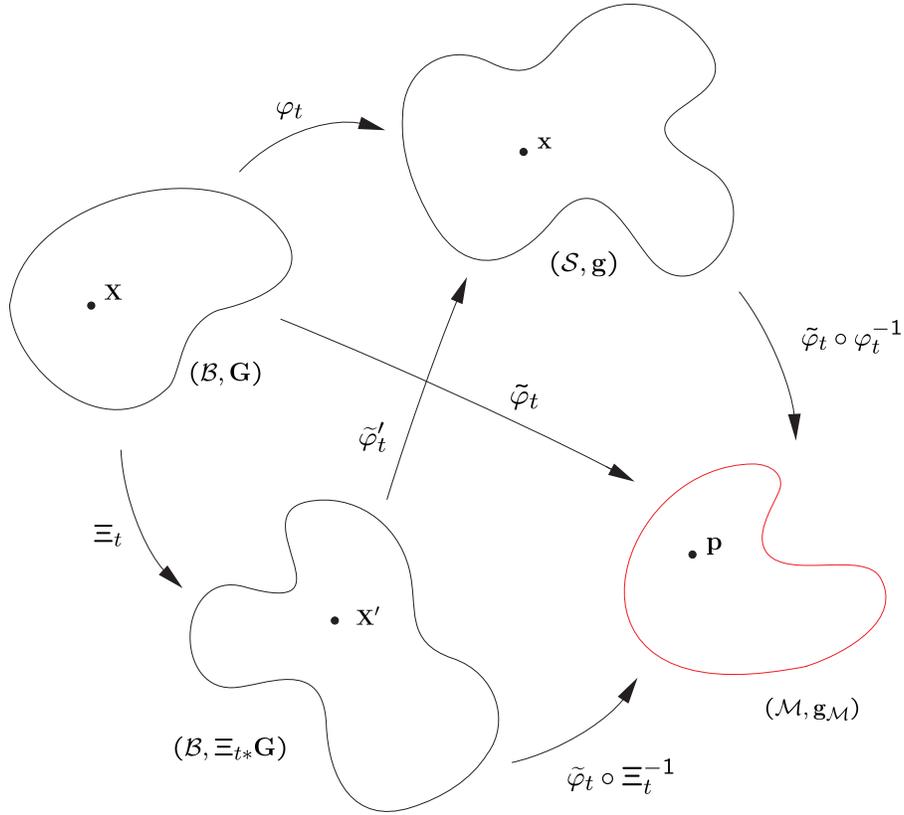}
\end{center}
\caption{\footnotesize Referential change of frame in a continuum
with microstructure.} \label{Cosserat_MaterialFraming}
\end{figure}
\paragraph{Change of Reference Frame.} A material change of frame is a
diffeomorphism
\begin{equation}
    \Xi_t:(\mathfrak{B},\mathbf{G})\rightarrow(\mathfrak{B},\mathbf{G}').
\end{equation}
A change of frame can be thought of as a change of coordinates in
the reference configuration (passive definition) or a
rearrangement of microstructure (active definition). Under such a
framing, a nice subset $\mathcal{U}$ is mapped to another nice
subset $\mathcal{U}'=\Xi_t(\mathcal{U})$ and a material point
$\mathbf{X}$ is mapped to $\mathbf{X}'=\Xi_t(\mathbf{X})$ (see
Fig. \ref{Cosserat_MaterialFraming}). The deformation mappings for
the new reference configuration are
$\varphi'_t=\varphi_t\circ\Xi_t^{-1}$ and
$\widetilde{\varphi}'_t=\widetilde{\varphi}_t\circ\Xi_t^{-1}$.
This can be clearly seen in Fig. \ref{Cosserat_MaterialFraming}.
The material velocity in $\mathcal{U}'$ is
\begin{equation}
    \mathbf{V}'(\mathbf{X}',t)=\frac{\partial}{\partial
    t}\varphi'_t(\mathbf{X}')=\frac{\partial \varphi_t}{\partial
    t}\circ\Xi_t^{-1}(\mathbf{X}')+T\varphi_t\circ\frac{\partial \Xi_t^{-1}}{\partial
    t}(\mathbf{X}'),
\end{equation}
where partial derivatives are calculated for fixed $\mathbf{X}'$.
We assume that
\begin{equation}
    \Xi_t\big|_{t=t_0}=id,~~~\frac{\partial \Xi_t}{\partial
    t}(\mathbf{X})=\mathbf{W}(\mathbf{X},t).
\end{equation}
Note that $\mathbf{W}$ is the infinitesimal generator of the
rearrangement $\Xi_t$. It is an easy exercise to show that
\begin{equation}
    \mathbf{V}'=\mathbf{V}\circ\Xi_t^{-1}-\mathbf{F}\mathbf{F}_{\Xi}^{-1}\cdot\mathbf{W}\circ\Xi_t^{-1}.
\end{equation}
Thus, at $t=t_0$
\begin{equation}
    \mathbf{V}'=\mathbf{V}-\mathbf{F}\mathbf{W}.
\end{equation}
Similarly
\begin{equation}
    \widetilde{\mathbf{V}}'=\widetilde{\mathbf{V}}-\widetilde{\mathbf{F}}\mathbf{W}.
\end{equation}
Note that
\begin{equation}
    \mathbf{G}' = (\varphi_t\circ\Xi_t^{-1})^*\circ\varphi_{t*}\mathbf{G}=
    (\Xi_t^{-1})^*\circ\varphi_t^*\circ\varphi_{t*}\mathbf{G}
   = (\Xi_t^{-1})^*\mathbf{G}=\Xi_{t*}\mathbf{G}=\left(T\Xi_t\right)^{-*}\mathbf{G}\left(T\Xi_t\right)^{-1}.
\end{equation}
And
\begin{equation}
    \mathbf{F}'=\Xi_{t*}\mathbf{F}=\mathbf{F}\circ (T\Xi_t)^{-1}.
\end{equation}
The material internal energy density is assumed to transform
tonsorially, i.e.
\begin{equation}\label{Cartan-Material}
    E'(t,\mathbf{X}',\mathbf{G}')=E(t,\mathbf{X},\mathbf{G}).
\end{equation}
This means that internal energy density at $\mathbf{X}'$ evaluated
by the transformed metric $\mathbf{G}'$ is equal to the internal
energy density at $\mathbf{X}$ evaluated by the metric
$\mathbf{G}$. We know that $\mathbf{G}'=\Xi_{t*}\mathbf{G}$, and
thus
\begin{equation}
    E'(t,\mathbf{X}',\mathbf{G})=E(t,\mathbf{X},\Xi_t^*\mathbf{G}).
\end{equation}
Therefore
\begin{equation}
    \frac{d}{dt}\Big|_{t=t_0}E'(t,\mathbf{X}',\mathbf{G})=\frac{\partial E}{\partial
    t}+\frac{\partial E}{\partial
    \mathbf{G}}:\mathfrak{L}_{\mathbf{W}}\mathbf{G}.
\end{equation}

\paragraph{Balance of Energy for Reframings of the Reference Configuration.} Consider a deformation mapping
$\varphi_t:\mathcal{B}\rightarrow\mathcal{S}$ and a referential
diffeomorphism $\Xi_t:\mathfrak{B}\rightarrow\mathfrak{B}$. The
mappings
$\varphi'_t=\varphi_t\circ\Xi_t^{-1}:\mathcal{B}'\rightarrow\mathcal{S}$
and
$\widetilde{\varphi}'_t=\widetilde{\varphi}_t\circ\Xi_t^{-1}:\mathcal{B}'\rightarrow\mathcal{M}$,
where $\mathcal{B}'=\Xi_t(\mathcal{B})$, represent the deformation
of the new (evolved) reference configuration. Balance of energy
for $\Xi_t(\mathcal{U})$ should include the following two groups
of terms:
\begin{itemize}
    \item [i)] Looking at $(\varphi'_t,\widetilde{\varphi}'_t)$ as the deformation of
    $\mathcal{B}'$ in $\mathcal{S}\times \mathcal{M}$, one has the usual material
    energy balance for $\Xi_t(\mathcal{U})$. Transformation of
    fields from $(\mathfrak{B},\mathbf{G})$ to $(\mathfrak{B},\mathbf{G}')$ follows Cartan's
    space-time theory.
    \item [ii)] Nonstandard terms may appear to represent the
    energy associated with the material evolution.
\end{itemize}
We expect to see some new terms that are work-conjugate to
$\mathbf{W}_t=\frac{\partial}{\partial t}\Xi_t$. Let us denote the
volume and surface forces conjugate to $\mathbf{W}$ by
$\mathbf{B}_0$ and $\mathbf{T}_0$, respectively.

Instead of looking at spatial framings, let us fix the deformed
configuration and look at framings of the reference configuration.
We postulate that energy balance for each nice subset
$\mathcal{U}'$ has the following form
\begin{eqnarray}
  && \frac{d}{d t}\int_{\mathcal{U}'}\rho'_0\left(E'+\frac{1}{2}\left\langle \! \left\langle\mathbf{V}',\mathbf{V}'\right\rangle \! \right\rangle
  +\frac{1}{2}J'\left\langle \!\! \left\langle\widetilde{\mathbf{V}}',\widetilde{\mathbf{V}}'\right\rangle \!\! \right\rangle\right)dV'=
    \int_{\mathcal{U}'}\rho'_0\left(\left\langle \! \left\langle\mathbf{B}',\mathbf{V}'\right\rangle \! \right\rangle+
    \left\langle \!\! \left\langle \widetilde{\mathbf{B}}',\widetilde{\mathbf{V}}'\right\rangle \!\! \right\rangle
    +R'\right)dV'   \nonumber \\
  && \label{material-energy-new} ~~~~~+ \int_{\partial\mathcal{U}'}\left(\left\langle \! \left\langle\mathbf{T}',\mathbf{V}'\right\rangle \! \right\rangle+
   \left\langle \!\! \left\langle\widetilde{\mathbf{T}}',\widetilde{\mathbf{V}}'\right\rangle \!\! \right\rangle+ H'\right)dA'
   +  \int_{\mathcal{U}'}\left\langle \! \left\langle\mathbf{B}'_0,\mathbf{W}_t\right\rangle \! \right\rangle dV'+
    \int_{\partial\mathcal{U}'}\left\langle \! \left\langle\mathbf{T}'_0,\mathbf{W}_t\right\rangle \! \right\rangle
    dA',
\end{eqnarray}
where $\mathcal{U}'=\Xi_t(\mathcal{U})$ and $\mathbf{B}'_0$ and
$\mathbf{T}'_0$ are unknown vector fields at this point. Using
Cartan's spacetime theory, it is assumed that the primed
quantities have the following relation with the unprimed
quantities
\begin{eqnarray}
  && dV'=\Xi_{t^*}dV,~~R'(\mathbf{X}',t)=R(\mathbf{X},t),~~\rho'_0(\mathbf{X}',t)=\rho_0(\mathbf{X}), \nonumber \\
  && H'(\mathbf{X}',\hat{\mathbf{N}}',t)=H(\mathbf{X},\hat{\mathbf{N}},t),~~J'=J,\\
  &&
  \mathbf{T}'(\mathbf{X}',\hat{\mathbf{N}}',t)=\mathbf{T}(\mathbf{X},\hat{\mathbf{N}},t),
  ~\widetilde{\mathbf{T}}'(\mathbf{X}',\hat{\mathbf{N}}',t)=\widetilde{\mathbf{T}}(\mathbf{X},\hat{\mathbf{N}},t).
  \nonumber
\end{eqnarray}
We assume that body force is transformed in such a way that
\begin{equation}
    \mathbf{B}'-\mathbf{A}'=\Xi_{t^*}(\mathbf{B}-\mathbf{A}),~~~\widetilde{\mathbf{B}}'-\widetilde{\mathbf{A}}'=\Xi_{t^*}(\widetilde{\mathbf{B}}-\widetilde{\mathbf{A}}).
\end{equation}
Thus
\begin{equation}
    (\mathbf{B}'-\mathbf{A}')\big|_{t=t_0}=\mathbf{B}-\mathbf{A},~~~(\widetilde{\mathbf{B}}'-\widetilde{\mathbf{A}}')\big|_{t=t_0}=\widetilde{\mathbf{B}}-\widetilde{\mathbf{A}}.
\end{equation}

Note that if $\alpha$ is a $3$-form on $\mathcal{U}$, then
\begin{equation}
    \frac{d}{dt}\Big|_{t=t_0}\int_{\mathcal{U}'}\alpha'=\int_{\mathcal{U}}\frac{d}{d t}\Big|_{t=t_0}\left(\Xi_t^*\alpha'\right),
\end{equation}
where $\mathcal{U}'=\Xi_t(\mathcal{U})$. Thus
\begin{equation}
    \frac{d}{dt}\Big|_{t=t_0}\int_{\mathcal{U}'}E'dV'=\int_{\mathcal{U}}\frac{d}{d t}\Big|_{t=t_0}\left(\Xi_t^*E'\right)dV
    =\int_{\mathcal{U}}\left(\frac{\partial E}{\partial t}+\frac{\partial E}{\partial
    \mathbf{G}}:\mathfrak{L}_{\mathbf{W}}\mathbf{G}\right)dV.
\end{equation}
Material energy balance for $\mathcal{U}'\subset\mathcal{B}'$ at
$t=t_0$ reads
\begin{eqnarray}
  && \int_{\mathcal{U}}\frac{\partial \rho_0}{\partial
t}\left(E+\frac{1}{2}\left\langle \!
\left\langle\mathbf{V}-\mathbf{F}\mathbf{W},
    \mathbf{V}-\mathbf{F}\mathbf{W}\right\rangle \! \right\rangle + \frac{1}{2}J \left\langle \!\!\left\langle\widetilde{\mathbf{V}}-\widetilde{\mathbf{F}}\mathbf{W},
    \widetilde{\mathbf{V}}-\widetilde{\mathbf{F}}\mathbf{W}\right\rangle \!\! \right\rangle \right)dV  \nonumber \\
  && ~~~~~+ \int_{\mathcal{U}}\rho_0\Bigg(\frac{\partial E}{\partial
    t}+\frac{\partial E}{\partial \mathbf{G}}:\mathfrak{L}_{\mathbf{W}}\mathbf{G}+\left\langle \!\! \left\langle\mathbf{V}-\mathbf{F}\mathbf{W},
    \mathbf{A}'\big|_{t=t_0}\right\rangle \!\! \right\rangle + J\left\langle \!\! \left\langle\widetilde{\mathbf{V}}-\widetilde{\mathbf{F}}\mathbf{W},
    \widetilde{\mathbf{A}}'\big|_{t=t_0}\right\rangle \!\! \right\rangle  \nonumber \\
  && ~~~~~~~~~~~~+\frac{1}{2}\frac{\partial J}{\partial t} \left\langle \!\!\left\langle\widetilde{\mathbf{V}}-\widetilde{\mathbf{F}}\mathbf{W},
    \widetilde{\mathbf{V}}-\widetilde{\mathbf{F}}\mathbf{W}\right\rangle \!\!   \right\rangle \Bigg)dV
    = \int_{\mathcal{U}}\rho_0\left(\left\langle \!\! \left\langle\mathbf{B}'\big|_{t=t_0},\mathbf{V}-\mathbf{F}\mathbf{W}\right\rangle \!\! \right\rangle+R\right)dV \nonumber \\
  && ~~~~~~+ \int_{\mathcal{U}}\rho_0 \left\langle \!\! \left\langle \widetilde{\mathbf{B}}'\big|_{t=t_0},\widetilde{\mathbf{V}}
  -\widetilde{\mathbf{F}}\mathbf{W}\right\rangle \!\! \right\rangle dV+
    \int_{\partial\mathcal{U}}\left(\left\langle \! \left\langle\mathbf{T},\mathbf{V}-\mathbf{F}\mathbf{W}\right\rangle \! \right\rangle+H\right)dA  \nonumber \\
  && ~~~~~~+\int_{\partial\mathcal{U}} \left\langle \!\! \left\langle\widetilde{\mathbf{T}},\widetilde{\mathbf{V}}-\widetilde{\mathbf{F}}\mathbf{W}\right\rangle \!\! \right\rangle dA
  +\int_{\mathcal{U}}\left\langle \! \left\langle\mathbf{B}_0,\mathbf{W}\right\rangle \! \right\rangle dV
  +\int_{\partial\mathcal{U}}\left\langle \! \left\langle\mathbf{T}_0,\mathbf{W}\right\rangle \! \right\rangle
  dA.
\end{eqnarray}
We know that $\mathbf{T}_0$ and $\mathbf{B}_0$ are defined on
$\mathcal{B}$ and $\mathbf{T}'_0$ and $\mathbf{B}'_0$ are the
corresponding quantities defined on $\Xi_t(\mathcal{B})$. Here we
assume that
\begin{equation}\label{transformation-material}
    \mathbf{T}'_0=\Xi_{t*}\mathbf{T}_0~~~~~\textrm{and}~~~~~\mathbf{B}'_0=\Xi_{t*}\mathbf{B}_0.
\end{equation}
Subtracting balance of energy for $\mathcal{U}$ from this and
noting that
$\left(\mathbf{A}'-\mathbf{B}'\right)_{t=t_0}=\mathbf{A}-\mathbf{B}$
and
$\left(\widetilde{\mathbf{A}}'-\widetilde{\mathbf{B}}'\right)_{t=t_0}=\widetilde{\mathbf{A}}-\widetilde{\mathbf{B}}$
one obtains
\begin{eqnarray} \label{covariant-material-energy-1}
  && \int_{\mathcal{U}}\frac{\partial \rho_0}{\partial t}\left(-\left\langle \! \left\langle\mathbf{V},
    \mathbf{F}\mathbf{W}\right\rangle \! \right\rangle+\frac{1}{2}\left\langle \! \left\langle\mathbf{F}\mathbf{W},
    \mathbf{F}\mathbf{W}\right\rangle \! \right\rangle-J\left\langle \!\! \left\langle\widetilde{\mathbf{V}},
    \widetilde{\mathbf{F}}\mathbf{W}\right\rangle \!\! \right\rangle+\frac{1}{2}J\left\langle \!\! \left\langle\widetilde{\mathbf{F}}\mathbf{W},
    \widetilde{\mathbf{F}}\mathbf{W}\right\rangle \!\! \right\rangle\right)dV \nonumber \\
  && ~~~~~~+
    \int_{\mathcal{U}}\rho_0\Bigg[\frac{\partial E}{\partial \mathbf{G}}:\mathfrak{L}_{\mathbf{W}}\mathbf{G}
    -\left\langle \! \left\langle\mathbf{F}\mathbf{W},\mathbf{A}\right\rangle \! \right\rangle
    -\left\langle \!\! \left\langle\widetilde{\mathbf{F}}\mathbf{W},J\widetilde{\mathbf{A}}\right\rangle \!\!
    \right\rangle \nonumber \\
  && ~~~~~~~~~~~~~~~~~~
    +\frac{\partial J}{\partial t}\left(-\left\langle \!\! \left\langle\widetilde{\mathbf{V}},
    \widetilde{\mathbf{F}}\mathbf{W}\right\rangle \!\! \right\rangle+\frac{1}{2}\left\langle \!\! \left\langle\widetilde{\mathbf{F}}\mathbf{W},
    \widetilde{\mathbf{F}}\mathbf{W}\right\rangle \!\! \right\rangle\right)    \Bigg]dV \nonumber \\
  && ~~~~~~=-\int_{\mathcal{U}}\left\langle \! \left\langle\rho_0\mathbf{B},\mathbf{F}\mathbf{W}\right\rangle \! \right\rangle dV-
    \int_{\partial\mathcal{U}}\left\langle \! \left\langle\mathbf{T},\mathbf{F}\mathbf{W}\right\rangle \! \right\rangle
    dA-\int_{\mathcal{U}}\left\langle \!\! \left\langle\rho_0\widetilde{\mathbf{B}},\widetilde{\mathbf{F}}\mathbf{W}\right\rangle \!\! \right\rangle dV \nonumber \\
   && ~~~~~~~~~
    -
    \int_{\partial\mathcal{U}}\left\langle \!\! \left\langle\widetilde{\mathbf{T}},\widetilde{\mathbf{F}}\mathbf{W}\right\rangle \!\! \right\rangle
    dA +\int_{\mathcal{U}}\left\langle \! \left\langle\mathbf{B}_0,\mathbf{W}\right\rangle \! \right\rangle dV
    +\int_{\partial\mathcal{U}}\left\langle \! \left\langle\mathbf{T}_0,\mathbf{W}\right\rangle \! \right\rangle dA.
\end{eqnarray}
We know that
\begin{equation}\label{Piola-Stress}
    \left\langle \! \left\langle\mathbf{T},\mathbf{F}\mathbf{W}\right\rangle \! \right\rangle
    =\left\langle \!\! \left\langle \mathbf{F}\mathbf{W},\left\langle \!\! \left\langle \mathbf{P},\hat{\mathbf{N}}\right\rangle \!\! \right\rangle \right\rangle \!\! \right\rangle
    ,~~~\left\langle \!\! \left\langle\widetilde{\mathbf{T}},\widetilde{\mathbf{F}}\mathbf{W}\right\rangle \!\! \right\rangle
    =\left\langle \!\! \left\langle \widetilde{\mathbf{F}}\mathbf{W},\left\langle \!\! \left\langle \widetilde{\mathbf{P}},\hat{\mathbf{N}}\right\rangle \!\! \right\rangle \right\rangle \!\! \right\rangle,
\end{equation}
where $\mathbf{P}$ is the first Piola-Kirchhoff stress tensor.
Thus, substituting (\ref{Piola-Stress}) into
(\ref{covariant-material-energy-1}), Cauchy's theorem implies that
\begin{equation}
    \mathbf{T}_0=\left\langle \!\! \left\langle\mathbf{P}_0,\hat{\mathbf{N}}\right\rangle \!\!
    \right\rangle,
\end{equation}
for some second-order tensor $\mathbf{P}_0$.
The surface integrals in material energy balance have the
following transformations (see \citet{YaMaOr2006} for a proof.)
\begin{equation}
    \int_{\partial\mathcal{U}}\left\langle \! \left\langle\mathbf{F}^{\textsf{T}}\mathbf{T},\mathbf{W}\right\rangle \! \right\rangle
    dA =  \int_{\mathcal{U}}\operatorname{Div}\left\langle \! \left\langle\mathbf{F}^{\textsf{T}}\mathbf{P},\mathbf{W}\right\rangle \! \right\rangle
    dV =  \int_{\mathcal{U}}\left[\left\langle \! \left\langle\operatorname{Div} (\mathbf{F}^{\textsf{T}}\mathbf{P}),\mathbf{W}\right\rangle \! \right\rangle+
   \mathbf{F}^{\textsf{T}}\mathbf{P}:\mathbf{\Omega}+\mathbf{F}^{\textsf{T}}\mathbf{P}:
    \mathbf{K}\right]dV.
    \label{CauchyThm}
\end{equation}
And
\begin{equation}
    \int_{\partial\mathcal{U}}\left\langle \!\! \left\langle\widetilde{\mathbf{F}}^{\textsf{T}}\widetilde{\mathbf{T}},\mathbf{W}\right\rangle \!\! \right\rangle
    dA =  \int_{\mathcal{U}}\operatorname{Div}\left\langle \!\! \left\langle\widetilde{\mathbf{F}}^{\textsf{T}}\widetilde{\mathbf{P}},\mathbf{W}\right\rangle \!\! \right\rangle
    dV =  \int_{\mathcal{U}}\left[\left\langle \!\! \left\langle\operatorname{Div} (\widetilde{\mathbf{F}}^{\textsf{T}}\mathbf{P}),\mathbf{W}\right\rangle \!\! \right\rangle+
   \widetilde{\mathbf{F}}^{\textsf{T}}\widetilde{\mathbf{P}}:\mathbf{\Omega}+\mathbf{F}^{\textsf{T}}\mathbf{P}:
    \mathbf{K}\right]dV,
    \label{CauchyThmMicro}
\end{equation}
where
\begin{eqnarray}
  && \mathbf{\Omega}_{IJ}=\frac{1}{2}\big(G_{IK}W^K{}_{|J}-G_{JK}W^K{}_{|I}\big)=\frac{1}{2}\left(W_{I|J}-W_{J|I}\right), \\
  &&
  \mathbf{K}_{IJ}=\frac{1}{2}\left(G_{IK}W^K{}_{|J}+G_{JK}W^K{}_{|I}\right)=
  \frac{1}{2}\left(W_{I|J}+W_{J|I}\right),~\mathbf{K}=\frac{1}{2}\mathfrak{L}_{\mathbf{W}}\mathbf{G}.
\end{eqnarray}
Similarly
\begin{equation}
    \int_{\partial\mathcal{U}}\left\langle \! \left\langle\mathbf{T}_0,\mathbf{W}\right\rangle \! \right\rangle
    dA =  \int_{\mathcal{U}}\operatorname{Div}\left\langle \! \left\langle \mathbf{P}_0,\mathbf{W}\right\rangle \! \right\rangle
    dV =  \int_{\mathcal{U}}\left[\left\langle \! \left\langle\operatorname{Div}\mathbf{P}_0,\mathbf{W}\right\rangle \! \right\rangle+
   \mathbf{P}_0:\mathbf{\Omega}+\mathbf{P}_0:\mathbf{K}\right]dV.
   \label{CauchyThm2}
\end{equation}
At time $t=t_0$ the transformed balance of energy should be the
same as the balance of energy for $\mathcal{U}$. Thus, subtracting
the material balance of energy for $\mathcal{U}$ from the above
balance law and considering conservation of mass and
micro-inertia, one obtains
\begin{eqnarray}
  && \int_{\mathcal{U}}\rho_0\frac{\partial E}{\partial\mathbf{G}}:
  \mathfrak{L}_{\mathbf{W}}\mathbf{G}~dV+\int_{\mathcal{U}}\left\langle \!\! \left\langle\rho_0\mathbf{F}^{\textsf{T}}
  \left(\mathbf{B}-\mathbf{A}\right),\mathbf{W}\right\rangle \!\! \right\rangle dV
  +\int_{\mathcal{U}}\left\langle \! \left\langle\rho_0\widetilde{\mathbf{F}}^{\textsf{T}}
  \left(\widetilde{\mathbf{B}}-\widetilde{\mathbf{A}}\right),\mathbf{W}\right\rangle \! \right\rangle dV  \nonumber \\
  && ~~~~~~~~~~~~-\int_{\mathcal{U}}\left\langle \! \left\langle \rho_0\mathbf{B}_{0},\mathbf{W}\right\rangle \! \right\rangle dV
  +\int_{\partial\mathcal{U}}\left\langle \!\! \left\langle\mathbf{F}^{\textsf{T}}\mathbf{T}+\widetilde{\mathbf{F}}^{\textsf{T}}\widetilde{\mathbf{T}}-\mathbf{T}_{0}
  ,\mathbf{W}\right\rangle \!\! \right\rangle dA=0.
\end{eqnarray}
Therefore
\begin{eqnarray}
  && \int_{\mathcal{U}}\left(2\rho_0\frac{\partial E}{\partial\mathbf{G}}
  +\mathbf{F}^{\textsf{T}}\mathbf{P}+\widetilde{\mathbf{F}}^{\textsf{T}}\widetilde{\mathbf{P}}-\mathbf{P}_0\right):\frac{1}{2}
  \mathfrak{L}_{\mathbf{W}}\mathbf{G}~dV
  +\int_{\mathcal{U}}\left(\mathbf{F}^{\textsf{T}}\mathbf{P}+\widetilde{\mathbf{F}}^{\textsf{T}}\widetilde{\mathbf{P}}-\mathbf{P}_0\right):\mathbf{\Omega}~dV
  \nonumber\\
  && \label{transformation-final}~~+\int_{\mathcal{U}}\Big\langle \!\! \Big\langle \rho_0\mathbf{F}^{\textsf{T}}\left(\mathbf{B}-\mathbf{A}\right)
  +\rho_0\widetilde{\mathbf{F}}^{\textsf{T}}
  \left(\widetilde{\mathbf{B}}-\widetilde{\mathbf{A}}\right)-\mathbf{B}_{0}
  +\operatorname{Div}\left(\mathbf{F}^{\textsf{T}}\mathbf{P}+\widetilde{\mathbf{F}}^{\textsf{T}}\widetilde{\mathbf{P}}\right)-\operatorname{Div}\mathbf{P}_0,\mathbf{W}\Big\rangle \!\!
\Big\rangle dV=0.
\end{eqnarray}
Using balance of linear and micro-linear momenta,
(\ref{transformation-final}) is simplified to read
\begin{eqnarray} \label{transformed1}
  && \int_{\mathcal{U}}\left(2\rho_0\frac{\partial E}{\partial\mathbf{G}}
  +\mathbf{F}^{\textsf{T}}\mathbf{P}+\widetilde{\mathbf{F}}^{\textsf{T}}\widetilde{\mathbf{P}}-\mathbf{P}_0\right):\frac{1}{2}
  \mathfrak{L}_{\mathbf{W}}\mathbf{G}~dV
  +\int_{\mathcal{U}}\left(\mathbf{F}^{\textsf{T}}\mathbf{P}+\widetilde{\mathbf{F}}^{\textsf{T}}\widetilde{\mathbf{P}}-\mathbf{P}_0\right):\mathbf{\Omega}~dV
  \nonumber
  \\
  && ~+\int_{\mathcal{U}}\left\langle \!\! \left\langle \operatorname{Div}\left(\mathbf{F}^{\textsf{T}}\mathbf{P}+\widetilde{\mathbf{F}}^{\textsf{T}}\widetilde{\mathbf{P}}-\mathbf{P}_0\right)
  -\mathbf{F}^{\textsf{T}}\operatorname{Div}\mathbf{P}-\widetilde{\mathbf{F}}^{\textsf{T}}\operatorname{Div}\widetilde{\mathbf{P}}-\mathbf{B}_{0},\mathbf{W}\right\rangle \!\! \right\rangle
  dV=0.
\end{eqnarray}
Because $\mathcal{U}$ and $\mathbf{W}$ are arbitrary, one obtains
\begin{eqnarray}
  && \mathbf{P}_0=2\rho_0\frac{\partial E}{\partial\mathbf{G}}
  +\mathbf{F}^{\textsf{T}}\mathbf{P}+\widetilde{\mathbf{F}}^{\textsf{T}}\widetilde{\mathbf{P}}, \label{transformation-1}\\
  &&
  \left(\mathbf{F}^{\textsf{T}}\mathbf{P}+\widetilde{\mathbf{F}}^{\textsf{T}}\widetilde{\mathbf{P}}-\mathbf{P}_0\right)^{\textsf{T}}
  =\mathbf{F}^{\textsf{T}}\mathbf{P}+\widetilde{\mathbf{F}}^{\textsf{T}}\widetilde{\mathbf{P}}-\mathbf{P}_0, \label{transformation-2} \\
  && \mathbf{B}_{0}=\operatorname{Div}\left(\mathbf{F}^{\textsf{T}}\mathbf{P}+\widetilde{\mathbf{F}}^{\textsf{T}}\widetilde{\mathbf{P}}-\mathbf{P}_0\right)
  -\mathbf{F}^{\textsf{T}}\operatorname{Div}\mathbf{P}-\widetilde{\mathbf{F}}^{\textsf{T}}\operatorname{Div}\widetilde{\mathbf{P}}.
\end{eqnarray}
Note that (\ref{transformation-2}) is trivially satisfied after
having (\ref{transformation-1}). Thus, we have
\begin{eqnarray}
  && \mathbf{P}_0=2\rho_0\frac{\partial E}{\partial\mathbf{G}}
  +\mathbf{F}^{\textsf{T}}\mathbf{P}+\widetilde{\mathbf{F}}^{\textsf{T}}\widetilde{\mathbf{P}}, \\
  &&
  \mathbf{B}_{0}=\operatorname{Div}\left(\mathbf{F}^{\textsf{T}}\mathbf{P}+\widetilde{\mathbf{F}}^{\textsf{T}}\widetilde{\mathbf{P}}-\mathbf{P}_0\right)
  -\mathbf{F}^{\textsf{T}}\operatorname{Div}\mathbf{P}-\widetilde{\mathbf{F}}^{\textsf{T}}\operatorname{Div}\widetilde{\mathbf{P}}.
\end{eqnarray}

\paragraph{Remark.} Note that $\mathbf{B}_{0}$ and
$\mathbf{P}_{0}$ are material tensors and hence the transformation
(\ref{transformation-material}) makes sense.

\vskip 0.1 in \noindent In summary, we have proven the following
theorem.

{\theorem Under a referential diffeomorphism
$\Xi_t:\mathfrak{B}\rightarrow\mathfrak{B}$, and assuming that
material energy density transforms tensorially, i.e.
\begin{equation}
    E'(t,\mathbf{X}',\mathbf{G})=E(t,\mathbf{X},\Xi_t^{*}\mathbf{G}),
\end{equation}
material energy balance has the following transformation
\begin{eqnarray}
  && \frac{d}{d
    t}\int_{\Xi_t(\mathcal{U})}\rho'_0\left(E'+\frac{1}{2}\left\langle \! \left\langle\mathbf{V}',\mathbf{V}'\right\rangle \! \right\rangle
    +\frac{1}{2}J'\left\langle \!\! \left\langle\widetilde{\mathbf{V}}',\widetilde{\mathbf{V}}'\right\rangle \!\! \right\rangle\right)dV' =
    \int_{\Xi_t(\mathcal{U})}\rho'_0\left(\left\langle \! \left\langle\mathbf{B}',\mathbf{V}'\right\rangle \! \right\rangle+
    \left\langle \!\! \left\langle\widetilde{\mathbf{B}}',\widetilde{\mathbf{V}}'\right\rangle \!\! \right\rangle
    +R'\right)dV' \nonumber \\
  && \label{Energy-Balance-Transformation} ~~+ \int_{\partial\Xi_t(\mathcal{U})}\left(\left\langle \! \left\langle\mathbf{T}',\mathbf{V}'\right\rangle \! \right\rangle
    +\left\langle \!\! \left\langle\widetilde{\mathbf{T}}',\widetilde{\mathbf{V}}'\right\rangle \!\!
    \right\rangle+H'\right)dA'+  \int_{\Xi_t(\mathcal{U})}\left\langle \! \left\langle\mathbf{B}'_0,\mathbf{W}_t\right\rangle \! \right\rangle
    dV'+\int_{\partial\Xi_t(\mathcal{U})}\left\langle \! \left\langle\mathbf{T}'_0,\mathbf{W}_t\right\rangle \! \right\rangle
    dA',
\end{eqnarray}
where
\begin{eqnarray}
  \mathbf{T}'_0 &=& \Xi_{t*}\left[\left\langle \!\!\! \left\langle 2\rho_0\frac{\partial E}{\partial\mathbf{G}}
  +\mathbf{F}^{\textsf{T}}\mathbf{P}+\widetilde{\mathbf{F}}^{\textsf{T}}\widetilde{\mathbf{P}},\hat{\mathbf{N}} \right\rangle \!\!\! \right\rangle \right], \\
  \mathbf{B}'_0 &=& \Xi_{t*}\left[\operatorname{Div}\left(\mathbf{F}^{\textsf{T}}\mathbf{P}+\widetilde{\mathbf{F}}^{\textsf{T}}\widetilde{\mathbf{P}}-\mathbf{P}_0\right)
  -\mathbf{F}^{\textsf{T}}\operatorname{Div}\mathbf{P}-\widetilde{\mathbf{F}}^{\textsf{T}}\operatorname{Div}\widetilde{\mathbf{P}}\right],
\end{eqnarray}
and the other quantities are already defined. }

\paragraph{\textbf{Consequences of Assuming Invariance of Energy
Balance.}} Let us now study the consequences of assuming material
covariance of energy balance. Material energy balance is invariant
under material diffeomorphisms if and only if the following
relations hold between the nonstandard terms
\begin{eqnarray}
  && \label{codition1} \mathbf{P}_0=\mathbf{0}~~~~~\textrm{or}~~~~~2\rho_0\frac{\partial E}{\partial\mathbf{G}}
  =-\mathbf{F}^{\textsf{T}}\mathbf{P}-\widetilde{\mathbf{F}}^{\textsf{T}}\widetilde{\mathbf{P}}, \\
  && \label{codition2}
  \mathbf{B}_{0}=\mathbf{0}~~~~~\textrm{or}~~~~~\operatorname{Div}\left(\mathbf{F}^{\textsf{T}}\mathbf{P}+\widetilde{\mathbf{F}}^{\textsf{T}}\widetilde{\mathbf{P}}\right)=
  \mathbf{F}^{\textsf{T}}\operatorname{Div}\mathbf{P}+\widetilde{\mathbf{F}}^{\textsf{T}}\operatorname{Div}\widetilde{\mathbf{P}}.
\end{eqnarray}

\vskip 0.3 in \subsection {\textbf{Covariant Elasticity for a
Special Class of Structured Continua}}

In this subsection, we consider two special types of structured
continua in which microstructure manifold is linked to reference
and ambient space manifolds. In the first example, we assume that
for any $\mathbf{X}\in\mathcal{B}$, microstructure manifold is
$(T_{\mathbf{X}}\mathcal{B},\mathbf{G})$. For such a continuum,
directors are ``attached" to material points. We call this
continuum a \emph{referentially constrained structured} (RCS)
continuum. In the second example, we assume that in the deformed
configuration, microstructure manifold for
$\mathbf{x}=\varphi_t(\mathbf{X})$ is
$(T_{\mathbf{x}}\mathcal{S},\mathbf{g})$. We call such a continuum
a \emph{spatially constrained structured} (SCS) continuum. For RCS
continua we look at both referential and spatial covariance of
energy balance. This is a concrete example of what we earlier
called a structured continuum with free microstructue. For SCS
continua we look at spatial covariance of energy balance.

As was mentioned earlier, in most treatments of continua with
microstructure, one has two balances of linear momenta; one for
standard forces and one for microstructure forces, and one balance
of angular momentum, which has contributions from both standard
and micro-forces. In this subsection, we show that in a special
case when microstructure manifold is the tangent space of the
ambient space manifold, one can obtain all the balance laws
covariantly using a single balance of energy. Interestingly, there
will be two balances of linear momenta and one balance of angular
momentum. We will also see that there are different possibilities
for defining ``covariance" and depending on what one calls
``covariance", balance laws have different forms.


\paragraph {\textbf{Materially Constrained Structured Continua.}} Given
$\mathbf{X}\in\mathcal{B}$, and
$\mathcal{M}=T_{\mathbf{X}}\mathcal{B}$, director velocity is
defined as
\begin{equation}
    \widetilde{\mathbf{V}}=\frac{\partial \widetilde{\varphi}_t(\mathbf{X})}{\partial t}.
\end{equation}
For writing energy balance in $\mathcal{S}$ we need to
push-forward the director velocity. The spatial director velocity
is defined as
\begin{equation}
    \widetilde{\mathbf{v}}=\varphi_{t*}\widetilde{\mathbf{V}}=\mathbf{F}\widetilde{\mathbf{V}}.
\end{equation}
Micro-traction $\widetilde{\mathbf{T}}$ has the coordinate
representation
\begin{equation}
    \widetilde{\mathbf{T}}=\widetilde{T}^A\mathbf{E}_A.
\end{equation}
Internal energy density has the form
$e=e(t,\mathbf{x},\mathbf{p}\circ\varphi_t^{-1},\mathbf{g},\mathbf{G}\circ\varphi_t^{-1})$.
Spatial and microstructure diffeomorphisms act on macro and
micro-forces independently as was explained in Section
\ref{sec:cov}. The resulting governing equations are exactly
similar to those obtained previously and thus we leave the
details.

\paragraph {\textbf{Spatially Constrained Structured Continua.}} In the previous
section we assumed that the standard ambient space and the
microstructure manifolds are independent in the sense that they
can have independent changes of frame. It seems that this is not
the case for most materials with microstructure and this is
perhaps why one sees only one balance of angular momentum, e.g. in
liquid crystals \citep{Ericksen1961, Leslie1968}. Here, we present
an example of a structured continuum in which the microstructure
manifold is linked to the standard ambient space manifold. We
assume that for each $\mathbf{x}\in\mathcal{S}$, the director at
$\mathbf{x}$, i.e., $\mathbf{p}(\mathbf{x})$ is an element of
$T_{\mathbf{x}}\mathcal{S}$. In other words
\begin{equation}
    \mathcal{M}_{\mathbf{x}}=T_{\mathbf{x}}\mathcal{S}~~~~~~~~\forall~\mathbf{x}\in\varphi_t(\mathcal{B}),
\end{equation}
i.e., for each $\mathbf{x}$ microstructure manifold is
$T_{\mathbf{x}}\mathcal{S}$ and $\widetilde{\varphi}$ is a
time-dependent vector in $T_{\mathbf{x}}\mathcal{S}$. In the fiber
bundle representation schematically shown in Fig.
\ref{Cosserat_Bundle}, this means that microstructure bundle is
$T\mathcal{S}$, i.e. the tangent bundle of the ambient space
manifold.

Here we assume that the director field is a single vector field.
Generalization of the results to cases where the director is a
tensor field would be straightforward. The microstructure
deformation gradient has the following representation
\begin{equation}
    \widetilde{\mathbf{F}}=T\widetilde{\varphi}_t\circ\mathbf{F},~~~~~
    \widetilde{\mathbf{F}}:T_{\mathbf{x}}\mathcal{S}\rightarrow T_{\mathbf{p}(\mathbf{x})}T_{\mathbf{x}}\mathcal{S}.
\end{equation}
In components
\begin{equation}
    \widetilde{\mathbf{F}}=\widetilde{F}^{a}{}_{b}~\mathbf{e}_{a}\otimes\mathbf{e}^{b}.
\end{equation}
Microstructure velocity is defined as
\begin{equation}
    \widetilde{\mathbf{v}}(\mathbf{x},t)=\frac{\partial}{\partial t}\Big|_{\mathbf{X}}\widetilde{\varphi}_t(\mathbf{x}).
\end{equation}
In components
\begin{equation}
    \widetilde{v}^a=\frac{\partial p^a}{\partial t}+\frac{\partial p^a}{\partial x^b}v^b+\gamma^a_{bc}v^bp^c.
\end{equation}
Or
\begin{equation}
    \widetilde{\mathbf{v}}=\dot{\mathbf{p}}=\frac{\partial \mathbf{p}}{\partial t}+\boldsymbol{\nabla}_{\mathbf{v}}\mathbf{p}.
\end{equation}

Now let us consider a spatial change of frame, i.e.
$\xi_t:\mathcal{S}\rightarrow\mathcal{S}$. Note that
$\varphi'_t=\xi_t\circ\varphi_t$ and because
$\widetilde{\varphi}\in T_{\mathbf{x}}\mathcal{S}$ we have
\begin{equation}
    \widetilde{\varphi}'_t(\mathbf{x}')=T\xi_t\cdot\widetilde{\varphi}_t(\mathbf{x}).
\end{equation}
Microstructure velocity in the new frame is defined as
\begin{equation}
    \widetilde{\mathbf{v}}'=\frac{\partial \mathbf{p}'}{\partial t}+\boldsymbol{\nabla}_{\mathbf{v}'}\mathbf{p}'.
\end{equation}
Noting that $\mathbf{p}'=\xi_{t*}\mathbf{p}$ and
$\mathbf{v}'=\xi_{t*}\mathbf{v}+\mathbf{w}_t$, we obtain
\begin{equation}
    \widetilde{\mathbf{v}}'=\frac{\partial }{\partial
    t}\Big|_{\mathbf{x}'}\left(\xi_{t*}\mathbf{p}\right)
    +\xi_{t*}\left(\boldsymbol{\nabla}_{\mathbf{v}}\mathbf{p}\right)+\boldsymbol{\nabla}_{\mathbf{w}}\left(\xi_{t*}\mathbf{p}\right).
\end{equation}
Note that\footnote{This can be proved as follows.
\begin{equation}
    \frac{\partial }{\partial t}\Big|_{\mathbf{x}}\mathbf{p}'=\frac{\partial }{\partial t}\Big|_{\mathbf{x}'}\mathbf{p}'
    +\frac{\partial }{\partial t}\Big|_{\mathbf{x}}\left[p'^{\alpha}(\xi(\mathbf{x}))\mathbf{e}_{\alpha}(\xi(\mathbf{x}))\right]
    =\left(\frac{\partial p'^{\alpha}}{\partial \xi^{\beta}}+\gamma^{\alpha}_{\lambda\beta}p'^{\lambda}\right)w^{\beta}\mathbf{e}_{\alpha}
    =\boldsymbol{\nabla}_{\mathbf{w}}\mathbf{p}'.
\end{equation}
}
\begin{equation}
    \frac{\partial }{\partial
    t}\Big|_{\mathbf{x}'}\left(\xi_{t*}\mathbf{p}\right)
    =\frac{\partial }{\partial t}\Big|_{\mathbf{x}}\left(\xi_{t*}\mathbf{p}\right)-\boldsymbol{\nabla}_{\mathbf{w}}\left(\xi_{t*}\mathbf{p}\right).
\end{equation}
Thus
\begin{equation}
    \widetilde{\mathbf{v}}'=\frac{\partial }{\partial
    t}\Big|_{\mathbf{x}}\left(\xi_{t*}\mathbf{p}\right)
    +\xi_{t*}\left(\boldsymbol{\nabla}_{\mathbf{v}}\mathbf{p}\right).
\end{equation}
Note also that
\begin{equation}
    \frac{\partial }{\partial
    t}\Big|_{\mathbf{x}}\left(\xi_{t*}\mathbf{p}\right)=\xi_{t*}\left(\frac{\partial \mathbf{p}}{\partial t}\right)
    +\boldsymbol{\nabla}_{\xi_{t*}\mathbf{p}}\mathbf{w}.
\end{equation}
Therefore
\begin{equation}
    \widetilde{\mathbf{v}}'=\xi_{t*}\widetilde{\mathbf{v}}+\boldsymbol{\nabla}_{\xi_{t*}\mathbf{p}}\mathbf{w}.
\end{equation}
This means that at time $t=t_0$
\begin{equation}
    \widetilde{\mathbf{v}}'=\widetilde{\mathbf{v}}+\boldsymbol{\nabla}_{\mathbf{p}}\mathbf{w}.
\end{equation}
We assume that microstructure body forces transform such that
$\widetilde{\mathbf{a}}'-\widetilde{\mathbf{b}}'=\xi_{t*}(\widetilde{\mathbf{a}}-\widetilde{\mathbf{b}})$.

For this structured continuum we assume that, in addition to
metric, internal energy density explicitly depends on a connection
too, i.e.\footnote{Note that this is similar to Palatini's
formulation of general relativity \citep{Wald1984}, where both
metric and connection are assumed to be fields.}
\begin{equation}
    e=e(t,\mathbf{x},\mathbf{p},\mathbf{g},\boldsymbol{\nabla}).
\end{equation}
The connection $\boldsymbol{\nabla}$ is assumed to be metric
compatible, i.e $\boldsymbol{\nabla}\mathbf{g}=\mathbf{0}$ but not
necessarily torsion-free, i.e., $\boldsymbol{\nabla}$ is not
necessarily the Levi-Civita connection. Therefore, under a change
of frame we have the following transformation of internal energy
density
\begin{equation}
    e'(t,\mathbf{x}',\mathbf{p}',\mathbf{g},\boldsymbol{\nabla})=e(t,\mathbf{x},\mathbf{p},\xi_t^*\mathbf{g},\xi_t^*\boldsymbol{\nabla}).
\end{equation}
Thus, at $t=t_0$
\begin{equation}
    \dot{\overline{e'}}=\dot{e}+\frac{\partial e}{\partial\mathbf{g}}:\mathfrak{L}_{\mathbf{w}}\mathbf{g}
    +\frac{\partial e}{\partial\boldsymbol{\nabla}}:\mathfrak{L}_{\mathbf{w}}\boldsymbol{\nabla}.
\end{equation}
We know that for a given connection $\boldsymbol{\nabla}$
\citep{MaHu1983}
\begin{equation}
    \mathfrak{L}_{\mathbf{w}}\boldsymbol{\nabla}=\boldsymbol{\nabla}\boldsymbol{\nabla}\mathbf{w}+\boldsymbol{\mathcal{R}}\cdot\mathbf{w}.
\end{equation}
Or in coordinates
\begin{equation}
    \left(\mathfrak{L}_{\mathbf{w}}\boldsymbol{\nabla}\right)^a{}_{bc}=w^a{}_{b|c}+\mathcal{R}^a{}_{dbc}w^d,
\end{equation}
where $\boldsymbol{\mathcal{R}}$ is the curvature tensor of
$(\mathcal{S},\mathbf{g})$.

Balance of energy for $\varphi_t(\mathcal{U})\subset \mathcal{S}$
is written as
\begin{eqnarray}
  && \frac{d}{dt}\int_{\varphi_t(\mathcal{U})}\rho(\mathbf{x},t)\left[e(t,\mathbf{x},\mathbf{p},\mathbf{g},\boldsymbol{\nabla})
   +\frac{1}{2}\left\langle \! \left\langle \mathbf{v},\mathbf{v} \right\rangle \! \right\rangle
   +\frac{1}{2}j\left\langle \! \left\langle \widetilde{\mathbf{v}},\widetilde{\mathbf{v}} \right\rangle \! \right\rangle
   \right]
   \nonumber \\
  && ~~~~=\int_{\varphi_t(\mathcal{U})}\rho(\mathbf{x},t)\left(\left\langle \! \left\langle \mathbf{b},\mathbf{v} \right\rangle \! \right\rangle
   +\left\langle \!\! \left\langle \widetilde{\mathbf{b}},\widetilde{\mathbf{v}}\right\rangle \!\! \right\rangle
   +r\right)+\int_{\partial\varphi_t(\mathcal{U})}\left(\left\langle \! \left\langle \mathbf{t},\mathbf{v} \right\rangle \!
  \right\rangle+ \left\langle \!\! \left\langle \widetilde{\mathbf{t}},\widetilde{\mathbf{v}} \right\rangle
  \!\!
  \right\rangle+h\right)da.
\end{eqnarray}
Let us postulate that energy balance is invariant under arbitrary
spatial changes of frame
$\xi_t:\mathcal{S}\rightarrow\mathcal{S}$, i.e.
\begin{eqnarray}
  && \frac{d}{dt}\int_{\varphi'_t(\mathcal{U})}\rho'(\mathbf{x}',t)\left[e'(t,\mathbf{x}',\mathbf{p}',\mathbf{g},\boldsymbol{\nabla})
   +\frac{1}{2}\left\langle \! \left\langle \mathbf{v}',\mathbf{v}' \right\rangle \! \right\rangle
   +\frac{1}{2}j'\left\langle \! \left\langle \widetilde{\mathbf{v}}',\widetilde{\mathbf{v}}' \right\rangle \! \right\rangle \right]     \nonumber \\
  && ~~~
  =\int_{\varphi'_t(\mathcal{U})}\rho'(\mathbf{x}',t)\left(\left\langle \! \left\langle \mathbf{b}',\mathbf{v}' \right\rangle \! \right\rangle
   +\left\langle \!\! \left\langle \widetilde{\mathbf{b}}',\widetilde{\mathbf{v}}'\right\rangle \!\! \right\rangle
   +r'\right) +\int_{\partial\varphi'_t(\mathcal{U})}\left(\left\langle \! \left\langle \mathbf{t}',\mathbf{v}' \right\rangle \!
  \right\rangle+ \left\langle \!\! \left\langle \widetilde{\mathbf{t}}',\widetilde{\mathbf{v}}' \right\rangle
  \!\!
  \right\rangle+h'\right)da.
\end{eqnarray}
We know that
\begin{eqnarray}
  && e'(t,\mathbf{x}',\mathbf{p}',\mathbf{g},\boldsymbol{\nabla})=e(t,\mathbf{x},\mathbf{p},\xi_t^{*}\mathbf{g},\xi_t^{*}\boldsymbol{\nabla}),~r'=r,~h'=h, \\
  && \rho'(\mathbf{x}',t)=\xi_{t*}\rho(\mathbf{x},t),~\mathbf{v}'=\xi_{t*}\mathbf{v}+\mathbf{w},~\mathbf{b}'-\mathbf{a}'=\xi_{t*}(\mathbf{b}-\mathbf{a}), \\
  &&
  \mathbf{t}'=\xi_{t*}\mathbf{t},~\widetilde{\mathbf{t}}'=\xi_{t*}\widetilde{\mathbf{t}},~\widetilde{\mathbf{b}}'-\widetilde{\mathbf{a}}'=\xi_{t*}(\widetilde{\mathbf{b}}-\widetilde{\mathbf{a}}).
\end{eqnarray}
Subtracting balance of energy for $\varphi_t(\mathcal{U})$ from
that of $\varphi'_t(\mathcal{U})$ at $t=t_0$, we obtain
\begin{eqnarray}\label{Energy}
  && \int_{\varphi_t(\mathcal{U})}\Bigg[\mathbf{L}_{\mathbf{v}}\rho\left(\frac{1}{2}\left\langle \! \left\langle \mathbf{w},\mathbf{w} \right\rangle \! \right\rangle
  +\left\langle \! \left\langle \mathbf{v},\mathbf{w}\right\rangle \!
  \right\rangle\right)+\mathbf{L}_{\mathbf{v}}(\rho j)\left(
  \frac{1}{2}\left\langle \! \left\langle \boldsymbol{\nabla}\mathbf{w}\cdot\mathbf{p},\boldsymbol{\nabla}\mathbf{w}\cdot\mathbf{p} \right\rangle \! \right\rangle
  +\left\langle \! \left\langle \widetilde{\mathbf{v}},\boldsymbol{\nabla}\mathbf{w}\cdot\mathbf{p}\right\rangle \!
  \right\rangle \right)\nonumber \\
  && ~~~~+\rho\left(\frac{\partial e}{\partial \mathbf{g}}:\mathfrak{L}_{\mathbf{w}}\mathbf{g}
  +\frac{\partial e}{\partial \boldsymbol{\nabla}}:(\boldsymbol{\nabla}\boldsymbol{\nabla}\mathbf{w}+\boldsymbol{\mathcal{R}}\cdot\mathbf{w})
  +\left\langle \! \left\langle \mathbf{a},\mathbf{w} \right\rangle \! \right\rangle
  +j\left\langle \! \left\langle \widetilde{\mathbf{a}},\boldsymbol{\nabla}\mathbf{w}\cdot\mathbf{p} \right\rangle \! \right\rangle \right)\Bigg]
       \nonumber \\
  && ~=\int_{\varphi_t(\mathcal{U})}\rho\left\langle \! \left\langle \mathbf{b},\mathbf{w} \right\rangle \! \right\rangle+
  \int_{\partial\varphi_t(\mathcal{U})} \left\langle \! \left\langle \mathbf{t},\mathbf{w} \right\rangle \! \right\rangle da
  +\int_{\varphi_t(\mathcal{U})}\left\langle \!\! \left\langle \rho\widetilde{\mathbf{b}},\boldsymbol{\nabla}\mathbf{w}\cdot\mathbf{p} \right\rangle \!\! \right\rangle
  +\int_{\partial\varphi_t(\mathcal{U})} \left\langle \!\! \left\langle \widetilde{\mathbf{t}},\boldsymbol{\nabla}\mathbf{w}\cdot\mathbf{p} \right\rangle \!\! \right\rangle da.
\end{eqnarray}
Assuming that $\xi_t$ is such that
$\widetilde{\mathbf{v}}'\big|_{t=t_0}-\widetilde{\mathbf{v}}=\mathbf{0}$,
i.e., $\boldsymbol{\nabla}\mathbf{w}=\mathbf{0}$, Cauchy's theorem
applied to (\ref{Energy}) implies that there is a second-order
tensor $\boldsymbol{\sigma}$ such that $\mathbf{t}=\left\langle \!
\left\langle \boldsymbol{\sigma},\mathbf{n} \right\rangle \!
\right\rangle$. Now applying Cauchy's theorem to (\ref{Energy})
for an arbitrary $\xi_t$ implies the existence of another
second-order tensor $\widetilde{\boldsymbol{\sigma}}$ such that
$\widetilde{\mathbf{t}}=\left\langle \! \left\langle
\widetilde{\boldsymbol{\sigma}},\mathbf{n} \right\rangle \!
\right\rangle$.

\paragraph{Remark.} Microstructure manifold is the tangent space
of the ambient space manifold at every point. However,
microstructure is not related to the deformation mapping. This is
why, unlike the so-called second-grade materials
(see \cite{FriedGurtin2006}), two separate stress tensors exist.

\vskip 0.2 in \noindent As $\mathcal{U}$ and $\mathbf{w}$ are
arbitrary, and replacing $\rho$ by $\rho dv$ in (\ref{Energy}), we
conclude that
\begin{eqnarray}
  && \mathbf{L}_{\mathbf{v}}\rho=0, \\
  && \mathbf{L}_{\mathbf{v}}j=0.
\end{eqnarray}
Now let us simplify the last two integrals in (\ref{Energy}). The
volume integral is simplified to read
\begin{equation}
    \int_{\varphi_t(\mathcal{U})}\left\langle \!\! \left\langle \rho(\widetilde{\mathbf{b}}-j\widetilde{\mathbf{a}}),\boldsymbol{\nabla}\mathbf{w}\cdot\mathbf{p} \right\rangle \!\! \right\rangle dv
    =\int_{\varphi_t(\mathcal{U})}\rho(\widetilde{\mathbf{b}}-j\widetilde{\mathbf{a}})\otimes\mathbf{p}:\left(\frac{1}{2}\mathfrak{L}_{\mathbf{w}}\mathbf{g}+\boldsymbol{\omega}\right)
    dv.
\end{equation}
The surface integral is simplified as
\begin{eqnarray}
  && \int_{\partial\varphi_t(\mathcal{U})} \left\langle \!\! \left\langle \widetilde{\mathbf{t}},\boldsymbol{\nabla}\mathbf{w}\cdot\mathbf{p} \right\rangle \!\! \right\rangle da
  = \int_{\varphi_t(\mathcal{U})}\left(\widetilde{\sigma}^{ad}p^c w_{a|c}\right)_{|d}dv \nonumber \\
  && ~~~~~~~~~~=
  \int_{\varphi_t(\mathcal{U})}\left[\left(\operatorname{div}\widetilde{\boldsymbol{\sigma}}\right)\otimes\mathbf{p}
  +\widetilde{\boldsymbol{\sigma}}\cdot\nabla\mathbf{p}\right]:\left(\frac{1}{2}\mathfrak{L}_{\mathbf{w}}\mathbf{g}+\boldsymbol{\omega}\right)dv
  +\int_{\varphi_t(\mathcal{U})}
  \widetilde{\sigma}^{ad}p^cw_{a|c|d}~dv, \nonumber\\
  && ~~~~~~~~~~=
  \int_{\varphi_t(\mathcal{U})}\left[\left(\operatorname{div}\widetilde{\boldsymbol{\sigma}}\right)\otimes\mathbf{p}
  +\widetilde{\boldsymbol{\sigma}}\cdot\nabla\mathbf{p}\right]:\left(\frac{1}{2}\mathfrak{L}_{\mathbf{w}}\mathbf{g}+\boldsymbol{\omega}\right)dv
  +\int_{\varphi_t(\mathcal{U})} \widetilde{\boldsymbol{\sigma}}\otimes\widetilde{\mathbf{p}}:\boldsymbol{\nabla}\boldsymbol{\nabla}\mathbf{w}~dv. \nonumber\\
  &&
\end{eqnarray}
%
%
%
Thus
\begin{eqnarray}
  && \int_{\varphi_t(\mathcal{U})}\left(-2\rho\frac{\partial e}{\partial \mathbf{g}}+\boldsymbol{\sigma}+\left(\operatorname{div}\widetilde{\boldsymbol{\sigma}}\right)\otimes\mathbf{p}
  +\widetilde{\boldsymbol{\sigma}}\cdot\nabla\mathbf{p}+\rho(\widetilde{\mathbf{b}}-j\widetilde{\mathbf{a}})\otimes\mathbf{p}\right):\frac{1}{2}\mathfrak{L}_{\mathbf{w}}\mathbf{g}
  ~dv \nonumber \\
  && ~~+\int_{\varphi_t(\mathcal{U})}\left(-2\rho\frac{\partial e}{\partial \mathbf{g}}+\boldsymbol{\sigma}+\left(\operatorname{div}\widetilde{\boldsymbol{\sigma}}\right)\otimes\mathbf{p}
  +\widetilde{\boldsymbol{\sigma}}\cdot\nabla\mathbf{p}+\rho(\widetilde{\mathbf{b}}-j\widetilde{\mathbf{a}})\otimes\mathbf{p}\right):\boldsymbol{\omega}~dv \nonumber \\
  && ~~+\int_{\varphi_t(\mathcal{U})}\left\langle \!\!\! \left\langle -\rho\mathbf{a}+\rho\mathbf{b}+\operatorname{div}\boldsymbol{\sigma}
  -\rho\frac{\partial e}{\partial \boldsymbol{\nabla}}:\boldsymbol{\mathcal{R}},\mathbf{w} \right\rangle \!\!\! \right\rangle
  dv, \nonumber\\
  && ~~+\int_{\varphi_t(\mathcal{U})} \left(-\rho\frac{\partial e}{\partial \boldsymbol{\nabla}}
  +\widetilde{\boldsymbol{\sigma}}\otimes\widetilde{\mathbf{p}}\right):\boldsymbol{\nabla}\boldsymbol{\nabla}\mathbf{w}~dv=0.
\end{eqnarray}
Therefore, because $\mathcal{U}, \mathbf{w},$ and $\mathbf{z}$ are
arbitrary we finally obtain
\begin{eqnarray}
  && \label{Covariance_MS1} \mathbf{L}_{\mathbf{v}}\rho=0, \\
  && \label{Covariance_MS11} \mathbf{L}_{\mathbf{v}}j=0, \\
  && \label{Covariance_MS2} \operatorname{div}\boldsymbol{\sigma}+\rho\mathbf{b}=\rho\mathbf{a}+\rho\frac{\partial e}{\partial \boldsymbol{\nabla}}:\boldsymbol{\mathcal{R}}, \\
  && \label{Covariance_MS4} 2\rho\frac{\partial e}{\partial \mathbf{g}}=\boldsymbol{\sigma}+\left(\operatorname{div}\widetilde{\boldsymbol{\sigma}}\right)\otimes\mathbf{p}
  +\widetilde{\boldsymbol{\sigma}}\cdot\nabla\mathbf{p}+\rho(\widetilde{\mathbf{b}}-j\widetilde{\mathbf{a}})\otimes\mathbf{p}, \\
  && \label{Covariance_MS5} \left[\boldsymbol{\sigma}+\left(\operatorname{div}\widetilde{\boldsymbol{\sigma}}\right)\otimes\mathbf{p}+\widetilde{\boldsymbol{\sigma}}\cdot \nabla
  \mathbf{p}+\rho\widetilde{\mathbf{b}}\otimes\mathbf{p}\right]^{\textsf{T}}=\boldsymbol{\sigma}+\left(\operatorname{div}\widetilde{\boldsymbol{\sigma}}\right)\otimes\mathbf{p}
  +\widetilde{\boldsymbol{\sigma}}\cdot\nabla\mathbf{p}+\rho(\widetilde{\mathbf{b}}-j\widetilde{\mathbf{a}})\otimes\mathbf{p},
  \nonumber \\
  && \\
  && \label{Covariance_MS6} \rho\frac{\partial e}{\partial \boldsymbol{\nabla}}=\widetilde{\boldsymbol{\sigma}}\otimes\mathbf{p}.
\end{eqnarray}
In component form, (\ref{Covariance_MS4}) reads
\begin{equation}
    2\rho\frac{\partial e}{\partial g_{ab}}
    =\sigma^{ab}+\widetilde{\sigma}^{ac}{}_{|c}p^b+\widetilde{\sigma}^{ac}p^b{}_{|c}
    +\rho(\widetilde{b}^a-\widetilde{a}^a)p^b=\sigma^{ab}+\rho(\widetilde{b}^a-\widetilde{a}^a)p^b+\left(\widetilde{\sigma}^{ac}p^b\right)_{|c}.
\end{equation}
Note that combining (\ref{Covariance_MS2}) and
(\ref{Covariance_MS6}), one can write balance of linear momentum
as
\begin{equation}
    \operatorname{div}\boldsymbol{\sigma}+\rho\mathbf{b}=\rho\mathbf{a}
    +\left(\widetilde{\boldsymbol{\sigma}}\otimes\mathbf{p}\right):\boldsymbol{\mathcal{R}}.
\end{equation}
This means that both stress and micro-stress tensors contribute to
balance of linear momentum. It is seen that there is a single
balance of linear momentum, a single balance of angular momentum
both with contributions from forces and micro-forces, and two
Doyle-Ericksen formulas.

We should mention that \citet{Toupin1962, Toupin1964} showed that
for elastic materials for which energy depends on gradient of the
deformation gradient, i.e. the second derivative of deformation
mapping, balance of linear momentum and angular momentum are both
coupled for micro and macro forces. However, as was mentioned
earlier, here we are not considering second-grade materials.

%
%

\vskip 0.2 in \paragraph{\textbf{Generalized Covariance of Energy
Balance for Spatially Constrained Structured Continua.}}

In all the previous examples we observed that covariance under a
single spatial diffeomorphism cannot lead to a separate balance of
micro-linear momentum. Let us consider two diffeomorphisms
$\xi_t,\eta_t:\mathcal{S}\rightarrow\mathcal{S}$ such that both
are identity at $t=t_0$ and
\begin{equation}
    \mathbf{z}\neq\mathbf{w}
    ,~\boldsymbol{\nabla}\mathbf{z}\neq\boldsymbol{\nabla}\mathbf{w}
    ,~\boldsymbol{\nabla}\boldsymbol{\nabla}\mathbf{z}=\boldsymbol{\nabla}\boldsymbol{\nabla}\mathbf{w},
\end{equation}
where
\begin{equation}
    \mathbf{w}=\frac{\partial}{\partial t}\Big|_{t=t_0}\xi_t,~\mathbf{z}=\frac{\partial}{\partial t}\Big|_{t=t_0}\eta_t.
\end{equation}
We assume that under the simultaneous actions of these two
diffeomorphisms, $\eta_t$ acts on micro-quantities and $\xi_t$
acts on the remaining quantities (including metric and
connection). Thus, in the new frame
\begin{equation}
    \mathbf{p}'=\eta_{t*}\mathbf{p},~\widetilde{\mathbf{v}}'=\widetilde{\mathbf{v}}+\boldsymbol{\nabla}_{\mathbf{p}}\mathbf{z}
    ,~\widetilde{\mathbf{a}}'-\widetilde{\mathbf{b}}'=\eta_{t*}(\widetilde{\mathbf{a}}-\widetilde{\mathbf{b}}).
\end{equation}
We assume that energy balance is invariant under the simultaneous
actions of $\xi_t$ and $\eta_t$ and call this a \emph{generalized
covariance}. Therefore, generalized covariance implies that at
time $t=t_0$
\begin{eqnarray}
  && \int_{\varphi_t(\mathcal{U})}\Bigg[\mathbf{L}_{\mathbf{v}}\rho\left(\frac{1}{2}\left\langle \! \left\langle \mathbf{w},\mathbf{w} \right\rangle \! \right\rangle
  +\left\langle \! \left\langle \mathbf{v},\mathbf{w}\right\rangle \!
  \right\rangle\right)
  +\mathbf{L}_{\mathbf{v}}(\rho j)\left(\frac{1}{2}\left\langle \! \left\langle \boldsymbol{\nabla}\mathbf{z}\cdot\mathbf{p},\boldsymbol{\nabla}\mathbf{z}\cdot\mathbf{p} \right\rangle \! \right\rangle
  +\left\langle \! \left\langle \widetilde{\mathbf{v}},\boldsymbol{\nabla}\mathbf{z}\cdot\mathbf{p}\right\rangle \!
  \right\rangle \right)\nonumber \\
  && ~~~~~~~~~~+\rho\left(\frac{\partial e}{\partial \mathbf{g}}:\mathfrak{L}_{\mathbf{w}}\mathbf{g}
  +\frac{\partial e}{\partial \boldsymbol{\nabla}}:(\boldsymbol{\nabla}\boldsymbol{\nabla}\mathbf{w}+\boldsymbol{\mathcal{R}}\cdot\mathbf{w})
  +\left\langle \! \left\langle \mathbf{a},\mathbf{w} \right\rangle \! \right\rangle+\left\langle \! \left\langle j\widetilde{\mathbf{a}},\boldsymbol{\nabla}\mathbf{z}\cdot\mathbf{p} \right\rangle \! \right\rangle \right)\Bigg]
       \nonumber \\
  && ~=\int_{\varphi_t(\mathcal{U})}\rho\left\langle \! \left\langle \mathbf{b},\mathbf{w} \right\rangle \! \right\rangle+
  \int_{\partial\varphi_t(\mathcal{U})} \left\langle \! \left\langle \mathbf{t},\mathbf{w} \right\rangle \! \right\rangle da
  +\int_{\varphi_t(\mathcal{U})}\left\langle \!\! \left\langle \rho\widetilde{\mathbf{b}},\boldsymbol{\nabla}\mathbf{z}\cdot\mathbf{p} \right\rangle \!\! \right\rangle
  +\int_{\partial\varphi_t(\mathcal{U})} \left\langle \!\! \left\langle \widetilde{\mathbf{t}},\boldsymbol{\nabla}\mathbf{z}\cdot\mathbf{p} \right\rangle \!\! \right\rangle da.
\end{eqnarray}
Arbitrariness of $\mathbf{w}$ and $\mathbf{z}$ gives us
conservation of mass $\mathbf{L}_{\mathbf{v}}\rho=0$, conservation
of microstructure inertia $\mathbf{L}_{\mathbf{v}}j=0$, and the
existence of stress tensors $\boldsymbol{\sigma}$ and
$\widetilde{\boldsymbol{\sigma}}$. Thus
\begin{eqnarray}
  && \int_{\varphi_t(\mathcal{U})}\rho\left(\frac{\partial e}{\partial \mathbf{g}}:\mathfrak{L}_{\mathbf{w}}\mathbf{g}
  +\frac{\partial e}{\partial \boldsymbol{\nabla}}:(\boldsymbol{\nabla}\boldsymbol{\nabla}\mathbf{w}+\boldsymbol{\mathcal{R}}\cdot\mathbf{w})
  +\left\langle \! \left\langle \mathbf{a},\mathbf{w} \right\rangle\!\right\rangle+\left\langle \! \left\langle j\widetilde{\mathbf{a}},\boldsymbol{\nabla}\mathbf{z}\cdot\mathbf{p} \right\rangle \! \right\rangle \right)
       \nonumber \\
  && ~~~~=\int_{\varphi_t(\mathcal{U})}\rho\left\langle \! \left\langle \mathbf{b},\mathbf{w} \right\rangle \! \right\rangle+
  \int_{\partial\varphi_t(\mathcal{U})} \left\langle \! \left\langle \mathbf{t},\mathbf{w} \right\rangle \! \right\rangle da
  +\int_{\varphi_t(\mathcal{U})}\left\langle \!\! \left\langle \rho\widetilde{\mathbf{b}},\boldsymbol{\nabla}\mathbf{z}\cdot\mathbf{p} \right\rangle \!\! \right\rangle \nonumber \\
  && ~~~~~~+
  \int_{\varphi_t(\mathcal{U})}\left[\left(\operatorname{div}\widetilde{\boldsymbol{\sigma}}\right)\otimes\mathbf{p}
  +\widetilde{\boldsymbol{\sigma}}\cdot\nabla\mathbf{p}\right]:\boldsymbol{\nabla}\mathbf{z}~dv
  +\int_{\varphi_t(\mathcal{U})} \widetilde{\boldsymbol{\sigma}}\otimes\widetilde{\mathbf{p}}:\boldsymbol{\nabla}\boldsymbol{\nabla}\mathbf{z}~dv.
\end{eqnarray}
Arbitrariness of $\mathbf{z},\mathbf{w}$, and $\mathcal{U}$, and
noting that
$\boldsymbol{\nabla}\boldsymbol{\nabla}\mathbf{z}=\boldsymbol{\nabla}\boldsymbol{\nabla}\mathbf{w}$,
one obtains
\begin{eqnarray}
  && \label{GCovariance_MS1} \operatorname{div}\boldsymbol{\sigma}+\rho\mathbf{b}=\rho\mathbf{a}+\rho\frac{\partial e}{\partial \boldsymbol{\nabla}}:\boldsymbol{\mathcal{R}}, \\
  && \label{GCovariance_MS2} 2\rho\frac{\partial e}{\partial \mathbf{g}}=\boldsymbol{\sigma}, \\
  && \label{GCovariance_MS3} \boldsymbol{\sigma}^{\textsf{T}}=\boldsymbol{\sigma},  \nonumber \\
  && \label{GCovariance_MS4} \operatorname{div}(\widetilde{\boldsymbol{\sigma}}\otimes\mathbf{p})
  +\rho\widetilde{\boldsymbol{b}}\otimes\mathbf{p}=\rho\widetilde{\boldsymbol{a}}\otimes\mathbf{p}, \\
  && \label{GCovariance_MS5} \rho\frac{\partial e}{\partial \boldsymbol{\nabla}}=\widetilde{\boldsymbol{\sigma}}\otimes\mathbf{p}.
\end{eqnarray}
It is seen that generalized covariance gives a separate balance of
micro-linear momentum, i.e. Eq.(\ref{GCovariance_MS4}).

%

\vskip 0.3 in
\section {\textbf{Examples of Continua with Microstructure}}

In this section, we present two examples of continua with
microstructure and obtain their governing equations covariantly.
We first look at a theory of elastic solids with voids (see
\cite{NunziatoCowin1979}), which is a structured continuum with a
one-dimensional microstructure manifold. We show that
microstructure covariance in this case gives all the balance laws
and a scalar Doyle-Ericksen formula. We then geometrically study
the classical theory of mixtures (see \cite{Bowen1967,
BedfordDrumheller1983, GreenNaghdi1967, Sampaio1976,
Williams1973}) and obtain the governing equations covariantly.

\subsection {\textbf{A Geometric Theory of Elastic Solids with Distributed Voids}}

An elastic solid with distributed voids can be thought of as a
structured continuum with a scalar microstructure kinematical
variable, as in \cite{Car1989}; here, we follow
\citet{NunziatoCowin1979}. In addition to the standard deformation
mapping, it is assumed that mass density has the following
multiplicative decomposition
\begin{equation}
    \rho_0(\mathbf{X})=\overline{\rho}_0(\mathbf{X},t)\nu_0(\mathbf{X},t),
\end{equation}
where $\overline{\rho}_0$ is the density of the matrix material
and $\nu_0$ is the matrix volume fraction and $0<\nu_0\leq 1$.
Deformation is a pair of mappings
$(\varphi_t,\widetilde{\varphi}_t):\mathcal{B}\times\mathcal{B}\rightarrow\mathcal{S}\times\mathbb{R}$.
Material void velocity and void deformation gradient (a one-form
on $\mathcal{B}$) are defined as
\begin{equation}
    \widetilde{V}(\mathbf{X},t)=\frac{\partial\nu_0(\mathbf{X},t)}{\partial t},
    ~~~\widetilde{\mathbf{F}}(\mathbf{X},t)=\frac{\partial \nu_0(\mathbf{X},t)}{\partial
    \mathbf{X}}.
\end{equation}
Spatial void velocity is defined as
$\widetilde{v}=\widetilde{V}\circ\varphi^{-1}$. Internal energy
density at $\mathbf{x}\in\mathcal{S}$ has the following form
\begin{equation}
    e=e(t,\mathbf{x},\mathbf{g},\nu,T\nu),
\end{equation}
where $\nu=\nu_0\circ\varphi$ and hence
\begin{equation}
    (T\nu)_a=\frac{\partial \nu}{\partial x^a}=F^{-A}{}_a\frac{\partial \nu}{\partial X^A}.
\end{equation}
For a subset $\varphi_t(\mathcal{U})\subset \mathcal{S}$, balance
of energy reads
\begin{eqnarray}\label{energy-balance-void}
  && \frac{d}{dt}\int_{\varphi_t(\mathcal{U})}\rho(\mathbf{x},t)\left(e(t,\mathbf{x},\mathbf{g},\nu,T\nu)
   +\frac{1}{2}\left\langle \! \left\langle \mathbf{v},\mathbf{v} \right\rangle \! \right\rangle
   +\frac{1}{2}\kappa~\widetilde{v}^{~2} \right)dv
   \nonumber \\
  && ~~~~=\int_{\varphi_t(\mathcal{U})}\rho(\mathbf{x},t)\left(\left\langle \! \left\langle \mathbf{b},\mathbf{v} \right\rangle \! \right\rangle
   + \widetilde{b}~\widetilde{v} +r\right)dv +\int_{\partial\varphi_t(\mathcal{U})}\left(\left\langle \! \left\langle \mathbf{t},\mathbf{v} \right\rangle \!
  \right\rangle+ \widetilde{t}~\widetilde{v}+h\right)da,
\end{eqnarray}
where $\kappa=\kappa(\mathbf{x},t)$ is the so-called equilibrated
inertia \citep{NunziatoCowin1979}, and $\widetilde{b}$ and
$\widetilde{t}$ are the void body force and traction,
respectively, and both are scalars.

Let us first consider a time-dependent spatial change of frame
$\xi_t:\mathcal{S}\rightarrow\mathcal{S}$ such that at $t=t_0$,
$\xi_{t_0}=id$. Under this change of frame
$\nu'(\mathbf{x}')=\nu(\mathbf{x})$ and hence
\begin{equation}
    e'=e'(t,\mathbf{x}',\mathbf{g},\nu',T\nu')=e(t,\mathbf{x},\xi_t^*\mathbf{g},\nu,T\nu).
\end{equation}
Therefore, at $t=t_0$
\begin{equation}
    \dot{\overline{e'}}=\dot{e}+\frac{\partial e}{\partial
    \mathbf{g}}:\mathfrak{L}_{\mathbf{w}}\mathbf{g},
\end{equation}
where $\mathbf{w}=\frac{\partial}{\partial t}\xi_t\big|_{t=t_0}$.
Subtracting balance of energy for $\varphi_t(\mathcal{U})$ from
that of $\varphi'_t(\mathcal{U})$ at $t=t_0$, gives the existence
of Cauchy stress and the standard balance laws \citep{YaMaOr2006}.

Let us now consider a microstructure change of frame
$\eta_t:(0,1]\rightarrow (0,1]$ such that $\eta_t\big|_{t=t_0}=id$
and
\begin{equation}
    \frac{\partial \eta_t(\nu)}{\partial t}=z_t(\nu).
\end{equation}
Void velocity in the new frame has the following form
\begin{equation}
    \widetilde{v}~'=\frac{\partial}{\partial
    t}\eta_t\circ\nu=\eta_{t*}\widetilde{v}+z_t.
\end{equation}
Thus, at $t=t_0$,
$\widetilde{v}~'(\nu)=\widetilde{v}(\nu)+z(\nu)$. Under the void
change of frame, we have
\begin{equation}
    e'(t,\mathbf{x},\mathbf{g},\nu',T\nu')=e(t,\mathbf{x},\mathbf{g},\nu,T\eta_t\cdot T\nu).
\end{equation}
Note that
\begin{equation}
    \frac{d}{dt}(T\eta_t\cdot T\nu)=\frac{d}{dt}\left(\frac{\partial \eta_t}{\partial \nu}\right)\frac{\partial \nu}{\partial
    \mathbf{X}}+\frac{\partial \eta_t}{\partial \nu}\frac{\partial \widetilde{v}}{\partial \mathbf{X}}
    =\frac{\partial z}{\partial \nu}\frac{\partial \nu}{\partial
    \mathbf{X}}+\frac{\partial^2 \eta_t(\nu)}{\partial \nu^2}\frac{\partial \nu}{\partial
    \mathbf{X}}+\frac{\partial \eta_t}{\partial
    \nu}\frac{\partial \widetilde{v}}{\partial \mathbf{X}}.
\end{equation}
Thus, at $t=t_0$
\begin{equation}
    \frac{d}{dt}(T\eta_t\cdot T\nu)
    =z'(\nu)\frac{\partial \nu}{\partial
    \mathbf{X}}+\frac{\partial \eta_t}{\partial
    \nu}\frac{\partial \widetilde{v}}{\partial \mathbf{X}}.
\end{equation}
Therefore, at $t=t_0$
\begin{equation}
    \dot{\overline{e'}}=\dot{e}+\frac{\partial e}{\partial \nu_{,A}}\nu_{,A}z'(\nu).
\end{equation}

Balance of energy in the new void frame at $t=t_0$ reads
\begin{eqnarray}\label{energy-balance-void1}
  && \int_{\varphi_t(\mathcal{U})}\rho\left[\dot{e}+\frac{\partial e}{\partial \nu_{,A}}\nu_{,A}z'(\nu)
   +\left\langle \! \left\langle \mathbf{v},\mathbf{a} \right\rangle \! \right\rangle
   +\frac{1}{2}\dot{\kappa} (\widetilde{v}+z)^2+\kappa \widetilde{a}(\widetilde{v}+z) \right]dv
   \nonumber \\
  && ~~~~=\int_{\varphi_t(\mathcal{U})}\rho\left(\left\langle \! \left\langle \mathbf{b},\mathbf{v} \right\rangle \! \right\rangle
   + \widetilde{b}(\widetilde{v}+z) +r\right)dv +\int_{\partial\varphi_t(\mathcal{U})}\left(\left\langle \! \left\langle \mathbf{t},\mathbf{v} \right\rangle \!
  \right\rangle+ \widetilde{t}(\widetilde{v}+z)+h\right)da.
\end{eqnarray}
Subtracting (\ref{energy-balance-void}) from
(\ref{energy-balance-void1}), one obtains
\begin{equation}\label{cov-void}
    \int_{\varphi_t(\mathcal{U})}\rho\left[\frac{\partial e}{\partial \nu_{,A}}\nu_{,A}z'(\nu)
   +\frac{1}{2}\dot{\kappa} (2\widetilde{v}z+z^2)+\kappa \widetilde{a}z \right]dv
   =\int_{\varphi_t(\mathcal{U})}\rho\widetilde{b}z dv +\int_{\partial\varphi_t(\mathcal{U})} \widetilde{t}z da.
\end{equation}
Because $z$ and $\mathcal{U}$ are arbitrary, we conclude that
$\dot{\kappa}=0$, which is the balance of equilibrated inertia
\citep{GoodmanCowin1972}. Using Cauchy's theorem in the above
identity, we conclude that there exists a vector field
$\widetilde{\boldsymbol{\sigma}}$ (void Cauchy stress), such that
$\widetilde{t}=\widetilde{\sigma}^a\hat{n}_a$. Therefore, the
surface integral in (\ref{cov-void}) can be simplified to read
\begin{equation}
    \int_{\partial\varphi_t(\mathcal{U})} \widetilde{t}z da
    =\int_{\varphi_t(\mathcal{U})}\left[(\operatorname{div}\widetilde{\boldsymbol{\sigma}})z
    +F^{-A}{}_a\widetilde{\sigma}^a\nu_{,A}z'\right]dv.
\end{equation}
Now, (\ref{cov-void}) can be rewritten as
\begin{equation}
    \int_{\varphi_t(\mathcal{U})}\left(\rho\frac{\partial e}{\partial \nu_{,A}}\nu_{,A}
   -F^{-A}{}_a~\widetilde{\sigma}^a\nu_{,A} \right)z'(\nu) dv
   -\int_{\varphi_t(\mathcal{U})}\left(\operatorname{div}\widetilde{\boldsymbol{\sigma}}
    +\rho\widetilde{b}-\rho\widetilde{a}\right)z(\nu) dv=0.
\end{equation}
Because $z$ and $z'$ can be chosen independently and $\mathcal{U}$
is arbitrary, we conclude that
\begin{eqnarray}
  && \label{linear-momentum-void}\operatorname{div}\widetilde{\boldsymbol{\sigma}}+\rho\widetilde{b}=\rho\widetilde{a}, \\
  && \label{DE-void} F^{-A}{}_a~\widetilde{\sigma}^a\nu_{,A}=\rho\frac{\partial e}{\partial
  \nu_{,A}}\nu_{,A}.
\end{eqnarray}
Eq. (\ref{linear-momentum-void}) is balance of equilibrated linear
momentum \citep{NunziatoCowin1979} and Eq. (\ref{DE-void}) is a
scalar Doyle-Ericksen formula.

\vskip 0.3 in
\subsection {\textbf{A Geometric Theory of Mixtures}}

In mixture theory, one is given a finite number of bodies
(constituents) that can penetrate into one another with the
understanding that there is no self penetration within a given
constituent. Here, for the sake of simplicity, we ignore diffusion
as our goal is to demonstrate the power of covariance arguments in
deriving the balance laws. We assume that in our mixture
$\textsf{M}$ there are two constituents; generalization of our
results to the case of $N$ constituents is straightforward. We
denote the constituents by $\textsf{1}$ and $\textsf{2}$. We
should mention that recently \citet{Mariano2005} studied some
invariance/covariance ideas for mixtures. Our approach is slightly
different as will be explained in the sequel.

Each constituent is assumed to have its own reference manifold
$\left({}^{\textsf{i}}\mathcal{B},{}^{\textsf{i}}\mathbf{G}\right),~\textsf{i}=\textsf{1},\textsf{2}$.
Deformation of $\textsf{M}$ is defined by two deformation mappings
${}^{\textsf{i}}\varphi_t,~\textsf{i}=\textsf{1},\textsf{2}$ such
that
\begin{equation}
    {}^{\textsf{i}}\varphi_t:\left({}^{\textsf{i}}\mathcal{B},{}^{\textsf{i}}\mathbf{G}\right)\rightarrow
    \left(\mathcal{S},{}^{\textsf{i}}\mathbf{g}\right),~\textsf{i}=\textsf{1},\textsf{2},
\end{equation}
i.e., it is assumed that the ambient space manifold $\mathcal{S}$
is equipped with two different metrics ${}^{\textsf{1}}\mathbf{g}$
and ${}^{\textsf{2}}\mathbf{g}$.\footnote{This is similar to what
\citet{Mariano2005} does when postulating covariance of energy
balance.} Material and spatial velocities are defined as
\begin{equation}
    {}^{\textsf{i}}\mathbf{V}(\mathbf{X}_{\textsf{i}},t)=\frac{\partial~ {}^{\textsf{i}}\varphi_t(\mathbf{X}_{\textsf{i}},t)}{\partial t}
    ,~~~{}^{\textsf{i}}\mathbf{v}={}^{\textsf{i}}\mathbf{V}\circ{}^{\textsf{i}}\varphi_t^{-1},~~~\textsf{i}=\textsf{1},\textsf{2}.
\end{equation}
Deformation gradients are tangent maps of the two deformation
mappings, i.e.,
${}^{\textsf{i}}\mathbf{F}=T~{}^{\textsf{i}}\varphi,~\textsf{i}=\textsf{1},\textsf{2}$.

Given
$\mathbf{x}\in{}^{\textsf{i}}\varphi_t(\mathcal{B}_{\textsf{i}})$,
it is assumed that this point is occupied by particles from both
$\mathcal{B}_{\textsf{1}}$ and $\mathcal{B}_{\textsf{2}}$, i.e.,
given a time $t_0$, $\mathbf{x}$ is the pre-image of particles
$\mathbf{X}_{\textsf{1}}$ and $\mathbf{X}_{\textsf{2}}$ defined as
\begin{equation}
    \mathbf{X}_{\textsf{1}}={}^{\textsf{1}}\varphi_{t_0}^{-1}(\mathbf{x})~~~\textrm{and}~~~\mathbf{X}_{\textsf{2}}={}^{\textsf{2}}\varphi_{t_0}^{-1}(\mathbf{x}).
\end{equation}
Thus, at a later time $t$
\begin{equation}
    {}^{\textsf{1}}\varphi_{t}(\mathbf{X}_{\textsf{1}})={}^{\textsf{1}}\varphi_{t}\circ{}^{\textsf{1}}\varphi_{t_0}^{-1}(\mathbf{x})\neq
    {}^{\textsf{2}}\varphi_{t}(\mathbf{X}_{\textsf{2}})={}^{\textsf{2}}\varphi_{t}\circ{}^{\textsf{2}}\varphi_{t_0}^{-1}(\mathbf{x}),
\end{equation}
i.e., in general, the two particles $\mathbf{X}_{\textsf{1}}$ and
$\mathbf{X}_{\textsf{2}}$ will occupy two different points of
$\mathcal{S}$ at time $t$. This means that one can have spatial
changes of frame that act separately on different
constituents.\footnote{This is closely related to what
\citet{Mariano2005} does in his energy balance covariance
argument.}

In the traditional formulation of mixture theories, for each
constituent, one assumes the existence of an internal energy
density and a ``growth" of internal energy density. Here, we
assume that each constituent has an internal energy that depends
on all the spatial metrics. For our two-phase mixture $\textsf{M}$
this means that
\begin{equation}
    e_{\textsf{1}}=e_{\textsf{1}}\left(t,\mathbf{x},{}^{\textsf{1}}\mathbf{g},{}^{\textsf{2}}\mathbf{g}\right)
    ~~~\textrm{and}~~~e_{\textsf{2}}=e_{\textsf{2}}\left(t,\mathbf{x},{}^{\textsf{1}}\mathbf{g},{}^{\textsf{2}}\mathbf{g}\right).
\end{equation}
Dependence of each internal energy density on both the spatial
metrics accounts for the interaction of constituents. Each
constituent is assumed to have its own mass density
$\rho_{\textsf{i}},~\textsf{i}=\textsf{1},\textsf{2}$ and mass
density at point $\mathbf{x}$ is defined as
\begin{equation}
    \rho(\mathbf{x},t)=\nu_{\textsf{1}}(\mathbf{x},t)\rho_{\textsf{1}}(\mathbf{x},t)+\nu_{\textsf{2}}(\mathbf{x},t)\rho_{\textsf{2}}(\mathbf{x},t),
\end{equation}
where $\nu_{\textsf{i}}$ are volume fractions of the constituents,
although at this point we do not need to define $\rho$.

Balance of energy for a subset
$\mathcal{U}_t={}^{\textsf{1}}\varphi_t(\mathcal{U}_{\textsf{1}})={}^{\textsf{2}}\varphi_t(\mathcal{U}_{\textsf{2}})\subset\mathcal{S}$
is written as
\begin{eqnarray}
  && \frac{d}{dt}\int_{\mathcal{U}_t}\sum_{\textsf{i}}\rho_{\textsf{i}}(\mathbf{x},t)\left[e_{\textsf{i}}\left(\mathbf{x},t,{}^{\textsf{1}}\mathbf{g},{}^{\textsf{2}}\mathbf{g}\right)
   +\frac{1}{2}\left\langle \! \left\langle {}^{\textsf{i}}\mathbf{v},{}^{\textsf{i}}\mathbf{v} \right\rangle \!  \right\rangle_{\textsf{i}}\right]
    \nonumber \\
  && ~~~=\int_{\mathcal{U}_t}\sum_{\textsf{i}}\rho_{\textsf{i}}(\mathbf{x},t)\left(\left\langle \! \left\langle {}^{\textsf{i}}\mathbf{b},{}^{\textsf{i}}\mathbf{v} \right\rangle \! \right\rangle_{\textsf{i}}
   +r_{\textsf{i}}\right)+\int_{\partial\mathcal{U}_t}\sum_{\textsf{i}}\left(\left\langle \! \left\langle {}^{\textsf{i}}\mathbf{t},{}^{\textsf{i}}\mathbf{v} \right\rangle \!
  \right\rangle_{\textsf{i}}+h_{\textsf{i}}\right)da,
\end{eqnarray}
where $\left\langle \! \left\langle .,. \right\rangle \!
\right\rangle_{\textsf{i}}$ is the inner product induced from the
metric ${}^{\textsf{i}}\mathbf{g}$ and all the other quantities
have the obvious meanings. Balance of energy can be simplified to
read
\begin{eqnarray}\label{mixture-energy-balance}
  && \int_{\mathcal{U}_t}\sum_{\textsf{i}}\mathbf{L}_{~{}^{\textsf{i}}\mathbf{v}}~\rho_{\textsf{i}}(\mathbf{x},t)\left[e_{\textsf{i}}\left(\mathbf{x},t,{}^{\textsf{1}}\mathbf{g}
  ,{}^{\textsf{2}}\mathbf{g}\right)
   +\frac{1}{2}\left\langle \! \left\langle {}^{\textsf{i}}\mathbf{v},{}^{\textsf{i}}\mathbf{v} \right\rangle \! \right\rangle_{\textsf{i}}\right]
   +\int_{\mathcal{U}_t}\sum_{\textsf{i}}\rho_{\textsf{i}}(\mathbf{x},t)\left[\dot{e}_{\textsf{i}}\left(\mathbf{x},t,{}^{\textsf{1}}\mathbf{g},{}^{\textsf{2}}\mathbf{g}\right)
   +\left\langle \! \left\langle {}^{\textsf{i}}\mathbf{v},{}^{\textsf{i}}\mathbf{a} \right\rangle \! \right\rangle_{\textsf{i}}\right]   \nonumber \\
  && ~~~~=\int_{\mathcal{U}_t}\sum_{\textsf{i}}\rho_{\textsf{i}}(\mathbf{x},t)\left(\left\langle \! \left\langle {}^{\textsf{i}}\mathbf{b},{}^{\textsf{i}}\mathbf{v} \right\rangle \! \right\rangle_{\textsf{i}}
   +r_{\textsf{i}}\right)+\int_{\partial\mathcal{U}_t}\sum_{\textsf{i}}\left(\left\langle \! \left\langle {}^{\textsf{i}}\mathbf{t},{}^{\textsf{i}}\mathbf{v} \right\rangle \!  \right\rangle_{\textsf{i}}+h_{\textsf{i}}\right)da.
\end{eqnarray}
Traditionally, a separate balance of energy is postulated for each
constituent \citep{Mariano2005}. Here, we only postulate a balance
of energy for the whole mixture.

We now consider a spatial diffeomorphism
$\xi_t:\mathcal{S}\rightarrow\mathcal{S}$ that acts only on
$(\mathcal{S},{}^{\textsf{1}}\mathbf{g})$ and is the identity map
at $t=t_0$. We postulate covariance of energy balance, i.e., in
the new spatial frame energy balance reads
\begin{eqnarray}\label{mixture-energy-balance-framing}
  && \frac{d}{dt}\int_{\mathcal{U}'_t}\rho'_{\textsf{1}}(\mathbf{x}',t)\left[e'_{\textsf{1}}\left(\mathbf{x}',t,{}^{\textsf{1}}\mathbf{g},{}^{\textsf{2}}\mathbf{g}\right)
   +\frac{1}{2}\left\langle \! \left\langle {}^{\textsf{1}}\mathbf{v}',{}^{\textsf{1}}\mathbf{v}' \right\rangle \!  \right\rangle_{\textsf{1}}\right]
    \nonumber \\
  && ~+\frac{d}{dt}\int_{\mathcal{U}'_t}\rho'_{\textsf{2}}(\mathbf{x}',t)\left[e'_{\textsf{2}}\left(\mathbf{x}',t,{}^{\textsf{1}}\mathbf{g},{}^{\textsf{2}}\mathbf{g}\right)
   +\frac{1}{2}\left\langle \! \left\langle {}^{\textsf{2}}\mathbf{v}',{}^{\textsf{2}}\mathbf{v}' \right\rangle \!  \right\rangle_{\textsf{2}}\right]
    \nonumber \\
  && ~~~=\int_{\mathcal{U}'_t}\rho'_{\textsf{1}}(\mathbf{x}',t)\left(\left\langle \! \left\langle {}^{\textsf{1}}\mathbf{b}',{}^{\textsf{1}}\mathbf{v}' \right\rangle \! \right\rangle_{\textsf{1}}
   +r'_{\textsf{1}}\right)+\int_{\partial\mathcal{U}'_t}\left(\left\langle \! \left\langle {}^{\textsf{1}}\mathbf{t}',{}^{\textsf{1}}\mathbf{v}' \right\rangle \!  \right\rangle_{\textsf{1}}+h'\right)da',  \nonumber \\
  && ~~~~+\int_{\mathcal{U}'_t}\rho'_{\textsf{2}}(\mathbf{x}',t)\left(\left\langle \! \left\langle {}^{\textsf{2}}\mathbf{b}',{}^{\textsf{2}}\mathbf{v}' \right\rangle \! \right\rangle_{\textsf{2}}
   +r'_{\textsf{2}}\right)+\int_{\partial\mathcal{U}'_t}\left(\left\langle \! \left\langle {}^{\textsf{2}}\mathbf{t}',{}^{\textsf{2}}\mathbf{v}' \right\rangle \!  \right\rangle_{\textsf{2}}+h'_{\textsf{2}}\right)da'.
\end{eqnarray}
Spatial velocities have the following transformations
\begin{equation}
    {}^{\textsf{1}}\mathbf{v}'=\xi_{t*}{}^{\textsf{1}}\mathbf{v}+\mathbf{w}_t~~~\textrm{and}~~~{}^{\textsf{2}}\mathbf{v}'={}^{\textsf{2}}\mathbf{v}.
\end{equation}
We assume that ${}^{\textsf{1}}\mathbf{b}$ is transformed such
that \citep{MaHu1983}
\begin{equation}
    {}^{\textsf{1}}\mathbf{b}'-{}^{\textsf{1}}\mathbf{a}'=\xi_{t*}\left({}^{\textsf{1}}\mathbf{b}-{}^{\textsf{1}}\mathbf{a}\right).
\end{equation}
Note also that
\begin{eqnarray}
  && e'_{\textsf{1}}\left(\mathbf{x}',t,{}^{\textsf{1}}\mathbf{g},{}^{\textsf{2}}\mathbf{g}\right)=e_{\textsf{1}}\left(\mathbf{x},t,\xi_t^*~{}^{\textsf{1}}\mathbf{g},{}^{\textsf{2}}\mathbf{g}\right), \\
  && e'_{\textsf{2}}\left(\mathbf{x}',t,{}^{\textsf{1}}\mathbf{g},{}^{\textsf{2}}\mathbf{g}\right)=e_{\textsf{2}}\left(\mathbf{x},t,\xi_t^*~{}^{\textsf{1}}\mathbf{g},{}^{\textsf{2}}\mathbf{g}\right).
\end{eqnarray}
Thus, at $t=t_0$
\begin{eqnarray}
  && \dot{\overline{e'_{\textsf{1}}}}=\dot{e}_{\textsf{1}}+\frac{\partial e_{\textsf{1}}}{\partial~{}^{\textsf{1}}\mathbf{g}}:\mathfrak{L}_{\mathbf{w}}{}^{\textsf{1}}\mathbf{g}, \\
  && \dot{\overline{e'_{\textsf{2}}}}=\dot{e}_{\textsf{2}}+\frac{\partial e_{\textsf{2}}}{\partial~{}^{\textsf{1}}\mathbf{g}}:\mathfrak{L}_{\mathbf{w}}{}^{\textsf{1}}\mathbf{g}.
\end{eqnarray}
Subtracting the energy balance (\ref{mixture-energy-balance}) from
(\ref{mixture-energy-balance-framing}) evaluated at $t=t_0$ yields
\begin{eqnarray}
  && \int_{\mathcal{U}_t}\mathbf{L}_{~{}^{\textsf{1}}\mathbf{v}}~\rho_{\textsf{1}}(\mathbf{x},t)\left[\left\langle \! \left\langle \mathbf{w},{}^{\textsf{1}}\mathbf{v} \right\rangle \!
   \right\rangle_{\textsf{1}}+\frac{1}{2}\left\langle \! \left\langle \mathbf{w},\mathbf{w} \right\rangle \!
   \right\rangle_{\textsf{1}}\right]
   +\int_{\mathcal{U}_t}\rho_{\textsf{1}}(\mathbf{x},t)\left(\frac{\partial e_{\textsf{1}}}{\partial~ {}^{\textsf{1}}\mathbf{g}}
   +\frac{\partial e_{\textsf{2}}}{\partial ~{}^{\textsf{1}}\mathbf{g}}\right):\mathfrak{L}_{\mathbf{w}}\mathbf{g}
    \nonumber \\
  && ~~~=\int_{\mathcal{U}_t}\rho_{\textsf{1}}(\mathbf{x},t)\left(\left\langle \! \left\langle {}^{\textsf{1}}\mathbf{b}-{}^{\textsf{1}}\mathbf{a},\mathbf{w} \right\rangle \! \right\rangle_{\textsf{1}}
   \right)+\int_{\partial\mathcal{U}_t}\left(\left\langle \! \left\langle {}^{\textsf{1}}\mathbf{t},\mathbf{w} \right\rangle \!
  \right\rangle_{\textsf{1}}\right)da.
\end{eqnarray}
Arbitrariness of $\mathcal{U}_t$ and $\mathbf{w}$ would guarantee
the existence of a Cauchy stress
${}^{\textsf{1}}\boldsymbol{\sigma}$ such that
${}^{\textsf{1}}\mathbf{t}=\langle\!\langle
{}^{\textsf{1}}\boldsymbol{\sigma},\hat{\mathbf{n}}
\rangle\!\rangle_{\textsf{1}}$ and also will give the following
after replacing $\rho_{\textsf{1}}$ by $\rho_{\textsf{1}}dv$
\begin{eqnarray}
  && \mathbf{L}_{~{}^{\textsf{1}}\mathbf{v}}~\rho_{\textsf{1}}=0,  \\
  && \operatorname{div}_{\textsf{1}}{}^{\textsf{1}}\boldsymbol{\sigma}+\rho_{\textsf{1}}{}^{\textsf{1}}\mathbf{b}=\rho_{\textsf{1}}{}^{\textsf{1}}\mathbf{a}, \\
  && {}^{\textsf{1}}\boldsymbol{\sigma}={}^{\textsf{1}}\boldsymbol{\sigma}^{\textsf{T}}, \\
  && {}^{\textsf{1}}\boldsymbol{\sigma}=2\rho_{\textsf{1}}\frac{\partial (e_{\textsf{1}}+e_{\textsf{2}})}{\partial~ {}^{\textsf{1}}\mathbf{g}}.
\end{eqnarray}
Similarly, assuming that $\xi_t:\mathcal{S}\rightarrow\mathcal{S}$
acts only on $(\mathcal{S},{}^{\textsf{2}}\mathbf{g})$ and
postulating energy balance covariance will give the following
balance laws.
\begin{eqnarray}
  && \mathbf{L}_{~{}^{\textsf{2}}\mathbf{v}}~\rho_{\textsf{2}}=0,  \\
  && \operatorname{div}_{\textsf{2}}{}^{\textsf{2}}\boldsymbol{\sigma}+\rho_{\textsf{2}}{}^{\textsf{2}}\mathbf{b}=\rho_{\textsf{2}}{}^{\textsf{2}}\mathbf{a}, \\
  && {}^{\textsf{2}}\boldsymbol{\sigma}={}^{\textsf{2}}\boldsymbol{\sigma}^{\textsf{T}}, \\
  && {}^{\textsf{2}}\boldsymbol{\sigma}=2\rho_{\textsf{2}}\frac{\partial (e_{\textsf{1}}+e_{\textsf{2}})}{\partial~ {}^{\textsf{2}}\mathbf{g}}.
\end{eqnarray}
Note the coupling in the Doyle-Ericksen formulas. Note also that
these balance laws can be pulled back to either
$\mathcal{B}_{\textsf{1}}$ or $\mathcal{B}_{\textsf{2}}$.

\vskip 0.3 in
\section {\textbf{Lagrangian Field Theory of Continua with Microstructure, Noether's Theorem and Covariance}}

The original formulations of Cosserat continua were mainly
variational \citep{Toupin1962, Toupin1964}. There have also been
recent geometric formulations in the literature \citep{CarMar2003,
FabMar2005}. In this section we consider a Lagrangian density that
depends explicitly on metrics and look at the corresponding
Euler-Lagrange equations. Then an explicit relation between
covariance and Noether's theorem is established. Similar to the
ambiguity encountered in covariant energy balance in terms of the
link of the microstructure manifold with the ambient space
manifold, here we will see that this ambiguity shows up in the
action of a given flow on different independent variables of the
Lagrangian density.

The Lagrangian may be regarded as a map
$L:T\mathcal{C}\rightarrow\mathbb{R}$, where $\mathcal{C}$ is the
space of some sections\footnote{See \cite{MaHu1983} for details in
the case of standard continua. The case of structured continua
would be a straightforward generalization.}, associated to the
Lagrangian density $\mathcal{L}$ and a volume element $dV(X)$ on
$\mathcal{B}$ and is defined as
\begin{equation}
    L(\varphi,\dot{\varphi},\widetilde{\varphi},\dot{\widetilde{\varphi}})
    =\int_{\mathcal{B}}\mathcal{L}\Big(X,\varphi(X),\dot{\varphi}(X),\mathbf{F}(X),\mathbf{G}(X),\mathbf{g}(\varphi(X))
    ,\widetilde{\varphi}(X),\dot{\widetilde{\varphi}}(X),\widetilde{\mathbf{F}}(X),\widetilde{\mathbf{g}}(\varphi(X))\Big)dV(X).
\end{equation}
Here $\varphi$ and $\widetilde{\varphi}$ are understood as fields
representing standard and microstructure deformations,
respectively. Note that, in general, one may need to consider more
than one microstructure field with possibly different tensorial
properties. Note also that in this material representation, the
two maps $\varphi$ and $\widetilde{\varphi}$ have the same role
and it is not clear from the Lagrangian density which one is the
standard deformation map. However, having the coordinate
representation for these two maps and their tangent maps, one can
see which one is the microstructure map. Note also that
$\mathbf{g}$ and $\widetilde{\mathbf{g}}$ are background metrics
with no dynamics.

The {\bfi action function} is defined as
\begin{equation}
    S(\varphi)=\int_{t_0}^{t_1}L(\varphi,\dot{\varphi},\widetilde{\varphi},\dot{\widetilde{\varphi}})dt.
\end{equation}
{\bfi Hamilton's principle} states that the physical configuration
$(\varphi,\widetilde{\varphi})$ is the critical point of the
action, i.e.
\begin{equation}
    \delta S=\mathbf{d}S(\varphi,\widetilde{\varphi})\cdot(\delta\varphi,\delta\widetilde{\varphi})=0.
\end{equation}
This can be simplified to read
\begin{eqnarray}
  && \int_{t_0}^{t_1}\int_{\mathcal{B}}\Big(\frac{\partial\mathcal{L}}{\partial\varphi}\cdot\delta\varphi+\frac{\partial\mathcal{L}}{\partial\dot{\varphi}}\cdot\delta\dot{\varphi}
    +\frac{\partial\mathcal{L}}{\partial\mathbf{F}}:\delta\mathbf{F}
    +\frac{\partial\mathcal{L}}{\partial\mathbf{g}}:\delta\mathbf{g}  \nonumber \\
  && ~~~~~~~~+\frac{\partial\mathcal{L}}{\partial\widetilde{\varphi}}\cdot\delta\widetilde{\varphi}+\frac{\partial\mathcal{L}}{\partial\dot{\widetilde{\varphi}}}\cdot\delta\dot{\widetilde{\varphi}}
    +\frac{\partial\mathcal{L}}{\partial\widetilde{\mathbf{F}}}:\delta\widetilde{\mathbf{F}}
    +\frac{\partial\mathcal{L}}{\partial\widetilde{\mathbf{g}}}:\delta\widetilde{\mathbf{g}} \Big)dV(X)dt=0.
\end{eqnarray}
As $\delta\varphi$ and $\delta\widetilde{\varphi}$ are
independent, we obtain the following Euler-Lagrange equations
\begin{eqnarray}
  && \frac{\partial\mathcal{L}}{\partial\varphi^a}-\frac{d}{d t}\frac{\partial\mathcal{L}}{\partial\dot{\varphi}^a}
    -\left(\frac{\partial\mathcal{L}}{\partial F^a{}_A}\right)_{|A}
    -\frac{\partial\mathcal{L}}{\partial F^b{}_A} F^c{}_A\gamma^b_{ac}+2\frac{\partial\mathcal{L}}{\partial g_{cd}}~g_{bd}\gamma^b_{ac}=0, \\
  && \nonumber \\
  && \frac{\partial\mathcal{L}}{\partial\widetilde{\varphi}^{\alpha}}-\frac{d}{d t}\frac{\partial\mathcal{L}}{\partial\dot{\widetilde{\varphi}^{\alpha}}}
    -\left(\frac{\partial\mathcal{L}}{\partial\widetilde{F}^{\alpha}{}_A}\right)_{|A}
    -\frac{\partial\mathcal{L}}{\partial\widetilde{F}^{\beta}{}_A} \widetilde{F}^{\mu}{}_A\widetilde{\gamma}^{\beta}_{\alpha\mu}
    +2\frac{\partial\mathcal{L}}{\partial \widetilde{g}_{\mu\lambda}}~\widetilde{g}_{\beta\lambda}\widetilde{\gamma}^{\beta}_{\alpha\mu}=0.
\end{eqnarray}
We know that because of material-frame-indifference, $\mathcal{L}$
depends on $\mathbf{F}$ and $\mathbf{g}$ through $\mathbf{C}$.
Thus, Euler-Lagrange equations for the standard deformation
mapping is simplified to read
\begin{eqnarray}
  && P_a{}^A{}_{|A}+\frac{\partial\mathcal{L}}{\partial\varphi^a}=\rho_0g_{ab}A^b, \\
  && \frac{\partial\mathcal{L}}{\partial\widetilde{\varphi}^{\alpha}}-\frac{d}{d t}\frac{\partial\mathcal{L}}{\partial\dot{\widetilde{\varphi}^{\alpha}}}
    +\widetilde{P}_{\alpha}{}^A{}_{|A}
    +\widetilde{P}_{\beta}{}^A \widetilde{F}^{\mu}{}_A\widetilde{\gamma}^{\beta}_{\alpha\mu}
    +2\frac{\partial\mathcal{L}}{\partial
    \widetilde{g}_{\mu\lambda}}~\widetilde{g}_{\beta\lambda}\widetilde{\gamma}^{\beta}_{\alpha\mu}=0,
\end{eqnarray}
where
\begin{equation}
    P_a{}^A=-\frac{\partial\mathcal{L}}{\partial F^a{}_A},~~~
    \widetilde{P}_{\alpha}{}^A=-\frac{\partial\mathcal{L}}{\partial\widetilde{F}^{\alpha}{}_A}.
\end{equation}
\vskip 0.2 in

When Euler-Lagrange equations are satisfied, given a symmetry of
the Lagrangian density \textbf{Noether's theorem} tells us what
its corresponding conserved quantity is. Suppose $\psi_s$ is a
flow on $\mathcal{S}$ generated by a vector field $\mathbf{w}$,
i.e.
\begin{equation}
    \frac{d}{ds}\Big|_{s=0}\psi_s\circ\varphi=\mathbf{w}\circ\varphi.
\end{equation}
Now if we assume that this flow leaves the microstructure
quantities unchanged, i.e., if we assume that the ambient space
manifold and the microstructure manifold are independent, then
invariance of the Lagrangian density means that
\begin{eqnarray}
  && \mathcal{L}\left(X^A,\psi_s^a(\varphi),\frac{\partial \psi_s^a}{\partial x^b}\dot{\varphi}^b,\frac{\partial \psi_s^a}{\partial x^b}F^b{}_A,G_{AB}
    ,-\frac{\partial \psi_s^c}{\partial x^a}\frac{\partial \psi_s^d}{\partial x^b}g_{cd},\widetilde{\varphi}^{\alpha},\dot{\widetilde{\varphi}}^{\alpha},\widetilde{F}^{\alpha}{}_A,\widetilde{g}_{\alpha\beta}\right)  \nonumber \\
  && ~~~~~~~~~~~~~~~~~~~~~~~~~~
  =\mathcal{L}\left(X^A,\varphi^a,\dot{\varphi}^a,F^a{}_A,G_{AB},g_{ab},\widetilde{\varphi}^{\alpha},\dot{\widetilde{\varphi}}^{\alpha},\widetilde{F}^{\alpha}{}_A,\widetilde{g}_{\alpha\beta}\right).
\end{eqnarray}
\citet{YaMaOr2006} proved that this implies the following two
conditions
\begin{eqnarray}
  && \label{DE} 2\frac{\partial\mathcal{L}}{\partial g_{ab}}=
    g^{bc}\frac{\partial\mathcal{L}}{\partial F^c{}_A} F^a{}_A
  +g^{bc}\frac{\partial\mathcal{L}}{\partial
  \dot{\varphi}^c}\dot{\varphi}^a, \\
  && \label{Spatial_Homogeneity}
  \frac{\partial\mathcal{L}}{\partial\varphi^a}=0,
\end{eqnarray}
i.e., the Doyle-Ericksen formula and spatial homogeneity of the
Lagrangian density.

\vskip 0.2 in \noindent Now, suppose $\eta_s$ is a flow on
$\mathcal{M}$ generated by a vector field $\mathbf{z}$, i.e.
\begin{equation}
    \frac{d}{ds}\Big|_{s=0}\eta_s\circ\widetilde{\varphi}=\mathbf{z}\circ\widetilde{\varphi}.
\end{equation}
Invariance of the Lagrangian density with respect to $\eta_s$
means that
\begin{eqnarray}
  && \mathcal{L}\left(X^A,\varphi^a,\dot{\varphi}^a,F^a{}_A,G_{AB},g_{ab},
    \eta_s^{\alpha}(\widetilde{\varphi}),\frac{\partial \eta_s^{\alpha}}{\partial p^{\beta}}~\dot{\widetilde{\varphi}}^{\beta},\frac{\partial \eta_s^{\alpha}}{\partial p^{\beta}}~\widetilde{F}^{\beta}{}_A,
    -\frac{\partial \eta_s^{\mu}}{\partial p^{\alpha}}\frac{\partial \eta_s^{\lambda}}{\partial p^{\beta}}~\widetilde{g}_{\mu\lambda}\right)  \nonumber \\
  && ~~~~~~~~~~~~~~~~~~~~~~~~~~
  =\mathcal{L}\left(X^A,\varphi^a,\dot{\varphi}^a,F^a{}_A,G_{AB},g_{ab},\widetilde{\varphi}^{\alpha},\dot{\widetilde{\varphi}}^{\alpha},\widetilde{F}^{\alpha}{}_A,\widetilde{g}_{\alpha\beta}\right).
\end{eqnarray}
Differentiating the above identity with respect to $s$ and
evaluating it for $s=0$, after some lengthy manipulations we
obtain
\begin{eqnarray}
  && \label{Invariance_1} 2\frac{\partial\mathcal{L}}{\partial \widetilde{g}_{\alpha\beta}}=
    \widetilde{F}^{\alpha}{}_A ~\widetilde{g}^{\beta\mu}\frac{\partial\mathcal{L}}{\partial\widetilde{F}^{\mu}{}_A}
  +g^{\beta\mu}\frac{\partial\mathcal{L}}{\partial
  \dot{\widetilde{\varphi}}^{\mu}}\dot{\widetilde{\varphi}}^{\alpha}, \\
  && \label{Invariance_2} \frac{d}{d t}\frac{\partial\mathcal{L}}{\partial \dot{\widetilde{\varphi}}^{\alpha}}
  +\left(\frac{\partial\mathcal{L}}{\partial\widetilde{F}^{\alpha}{}_A}\right)_{|A}-\frac{\partial\mathcal{L}}{\partial \dot{\widetilde{\varphi}}^{\lambda}}
  ~\widetilde{\gamma}^{\lambda}_{\alpha\mu}\dot{\widetilde{\varphi}}^{\mu}=0.
\end{eqnarray}
Assuming that $\mathcal{L}$ has the following splitting in terms
of internal energy density and kinetic energy
\begin{equation}
    \mathcal{L}=\rho_0e+\frac{1}{2}\rho_0\left\langle \! \left\langle \mathbf{V},\mathbf{V} \right\rangle \! \right\rangle_{\mathbf{g}}
    +\frac{1}{2}\widetilde{\rho}_0\left\langle \!\! \left\langle \widetilde{\mathbf{V}},\widetilde{\mathbf{V}} \right\rangle \!\!
    \right\rangle_{\widetilde{\mathbf{g}}},
\end{equation}
(\ref{Invariance_1}) is simplified to read
\begin{equation}
    2\rho_0\frac{\partial e}{\partial \widetilde{g}_{\alpha\beta}}=
    \widetilde{F}^{\alpha}{}_A ~\widetilde{g}^{\beta\mu}\frac{\partial\mathcal{L}}{\partial\widetilde{F}^{\mu}{}_A}
    =\widetilde{F}^{\alpha}{}_A \widetilde{P}^{\beta A}.
\end{equation}
Now let us simplify this relation and show that it is exactly
equivalent to (\ref{DE-M}). Note that
\begin{equation}
    \widetilde{P}^{\alpha A}=J\left(\mathbf{F}^{-1}\right)^A{}_b~\widetilde{\sigma}^{\alpha b}.
\end{equation}
Thus
\begin{equation}
    \widetilde{F}^{\alpha}{}_A \widetilde{P}^{\beta A}=J \left(\widetilde{\mathbf{F}}\mathbf{F}^{-1}\right)^{\alpha}{}_b~\widetilde{\sigma}^{\beta b}.
\end{equation}
Hence
\begin{equation}
    2\rho\frac{\partial e}{\partial \widetilde{g}_{\alpha\beta}}=\left(\mathbf{F}_0\right)^{\alpha}{}_b~\widetilde{\sigma}^{\beta b}.
\end{equation}
This means that (\ref{Invariance_1}) is equivalent to
(\ref{DE-M})!

Following \citet{YaMaOr2006}, it can be shown that using
Euler-Lagrange equations and some lengthy manipulations,
(\ref{Invariance_2}) can be simplified to read
\begin{equation}
    \frac{\partial \mathcal{L}}{\partial \widetilde{\varphi}^{\alpha}}=0.
\end{equation}
This means that if Lagrangian density is microstructurally
covariant, then it has to be microstructurally homogenous and a
micro-Doyle-Ericksen formula should be satisfied.

\citet{FabMar2005} study invariance of Lagrangian density of a
structured continuum under different groups of transformations. In
particular, they require invariance of the Lagrangian density when
the same copy of $SO(3)$ acts on ambient space and microstructure
manifolds in order to obtain balance of angular momentum. This
seems to be a matter of choice at first sight but can also be
understood as an interpretation of balance of angular momenta for
a special class of structured continua. In the following, we study
a similar symmetry of the Lagrangian density.

\paragraph{Constrained Microstructure Manifold.} Now let us assume
that $
\mathcal{M}(\mathbf{X})=T_{\varphi_t(\mathbf{X})}\mathcal{S}$. In
this case the Euler-Lagrange equations are
\begin{eqnarray}
  && P_a{}^A{}_{|A}+\frac{\partial\mathcal{L}}{\partial\varphi^a}=\rho_0g_{ab}A^b, \\
  && \nonumber \\
  && \frac{\partial\mathcal{L}}{\partial\widetilde{\varphi}^{a}}
  -\frac{d}{d t}\frac{\partial\mathcal{L}}{\partial\dot{\widetilde{\varphi}^a}}
    -\left(\frac{\partial\mathcal{L}}{\partial\widetilde{F}^a{}_A}\right)_{|A}
    -\frac{\partial\mathcal{L}}{\partial\widetilde{F}^b{}_A}\widetilde{F}^{c}{}_A\gamma^{b}_{ac}=0.
\end{eqnarray}
Now a flow on $\mathcal{S}$ would affect the microstructure
quantities too. In this case invariance of the Lagrangian density
means that
\begin{eqnarray}
  && \mathcal{L}\left(X^A,\psi_s^a(\varphi),\frac{\partial \psi_s^a}{\partial x^b}\dot{\varphi}^b,\frac{\partial \psi_s^a}{\partial x^b}F^b{}_A,G_{AB}
    ,-\frac{\partial \psi_s^c}{\partial x^a}\frac{\partial \psi_s^d}{\partial x^b}g_{cd},\frac{\partial \psi_s^a}{\partial x^b}\widetilde{\varphi}^{b}
    ,\frac{\partial \psi_s^a}{\partial x^b}\dot{\widetilde{\varphi}}^{b},\frac{\partial \psi_s^a}{\partial x^b}\widetilde{F}^{b}{}_A\right)  \nonumber \\
  && ~~~~~~~~~~~~~~~~~~~~~~~~~~
  =\mathcal{L}\left(X^A,\varphi^a,\dot{\varphi}^a,F^a{}_A,G_{AB},g_{ab},\widetilde{\varphi}^{a},\dot{\widetilde{\varphi}}^{a},\widetilde{F}^{a}{}_A\right).
\end{eqnarray}
Differentiating the above identity with respect to $s$ and
evaluating it at $s=0$ yields
\begin{eqnarray}
  && 2\frac{\partial\mathcal{L}}{\partial g_{ab}}=g^{bc}\left(\frac{\partial\mathcal{L}}{\partial F^c{}_A}F^a{}_A
     +\frac{\partial \mathcal{L}}{\partial \dot{\varphi}^c}\dot{\varphi}^a
     +\frac{\partial\mathcal{L}}{\partial \widetilde{\varphi}^c}\widetilde{\varphi}^a
     +\frac{\partial\mathcal{L}}{\partial \dot{\widetilde{\varphi}}^c}\dot{\widetilde{\varphi}}^a
     +\frac{\partial\mathcal{L}}{\partial \widetilde{F}^c{}_A}\widetilde{F}^a{}_A\right), \\
  && \frac{\partial \mathcal{L}}{\partial
  \varphi^a}-\left(\frac{\partial\mathcal{L}}{\partial \widetilde{\varphi}^c}\widetilde{\varphi}^a
     +\frac{\partial\mathcal{L}}{\partial \dot{\widetilde{\varphi}}^c}\dot{\widetilde{\varphi}}^a
     +\frac{\partial\mathcal{L}}{\partial \widetilde{F}^c{}_A}\widetilde{F}^a{}_A\right)\gamma^c_{ab}=0.
\end{eqnarray}


\vskip 0.4 in
\section {\textbf{Concluding Remarks}}

This paper first critically reviewed the geometry of structured
continua. Similar to classical continuum mechanics, one assumes
the existence of a well-defined reference configuration and each
material point  is mapped to its current position in the ambient
space by the standard deformation mapping. In addition to this,
each material point is given a director, which lies in a
microstructure manifold. A separate map, the microstructure
deformation mapping, maps each material point to its director,
which could be a scalar field, a vector field, or in general a
tensor field.

The Green-Naghdi-Rivlin Theorem relates balance laws to invariance
of balance of energy under some groups of transformations.
Previous attempts to extend this theorem to structured continua
were critically reviewed. It was explained that any generalization
of this theorem explicitly depends on the nature of the
microstructure manifold. It turns out that in most continua with
microstructure, the microstructure manifold is linked to the
ambient space manifold. We gave a concrete example of a structured
continuum, in which the ambient space is Euclidean, for which the
microstructure manifold is again $\mathbb{R}^3$ but thought of as
the tangent space of $\mathbb{R}^3$ at a given point. Postulating
balance of energy and its invariance under isometries of
$\mathbb{R}^3$, we obtained conservation of mass, balance of
linear momentum and balance of angular momentum with contributions
from both macro and micro-forces. Limiting oneself to rigid
motions does not allow one to obtain a separate balance of
micro-linear momentum. This leads one to think about investigating
covariant balance laws for structured continua.

We first assumed that the structured continuum is such that the
ambient space and the macrostructure manifold can have independent
reframings. We showed that postulating energy balance and its
invariance under spatial and microstructure diffeomorphisms gives
conservation of mass, existence of Cauchy stress and micro-Cauchy
stress, balance of linear and micro-linear momenta, balance of
angular and micro-angular momenta and two Doyle-Ericksen formulas.
We then considered structured continua for which the
microstructure manifold is somewhat constrained in the sense that
a spatial change of frame affects the microstructure quantities
too. As concrete examples, we defined materially and spatially
constrained structured continua. In a spatially constrained
continuum the microstructure bundle is the tangent bundle of the
ambient space manifold. In a materially constrained structured
continuum, microstructure manifold at a given point
$\mathbf{X}\in\mathcal{B}$ is $T_{\mathbf{X}}\mathcal{B}$. We
showed that postulating energy balance and its invariance under
spatial diffeomorphisms for a MCS continuum gives conservation of
mass, two balances of linear momentum, two balances of angular
momentum and two Doyle-Ericksen formulas. For a SCS continuum,
spatial covariance gives balances of linear and angular momenta,
which both have contributions from macro and micro forces. We then
defined a generalized covariance in which two separate maps act on
macro and micro quantities simultaneously. Under some assumptions,
we showed that generalized covariance can give a coupled balance
of angular momentum and two separate balances of linear momentum
for macro and micro forces.

As concrete examples of structured continua, we looked at elastic
solids with distributed voids and mixtures and obtained their
balance laws covariantly.

In the last part of the paper, we reviewed the Lagrangian field
theory of structured continua, when both ambient space and
microstructure manifolds are equipped with their own metrics.
Assuming that standard deformation mapping and microstructure
deformation mapping are independent, they would have independent
variations and hence Hamilton's principle of least action gives us
two sets of Euler-Lagrange equations. We then studied the
connection between Noether's theorem and covariance. It was
observed that there is some ambiguity in making this connection.
The ambiguity arises from the fact that there are different
possibilities in defining covariance for a Lagrangian density. One
choice is to assume that the Lagrangian density is covariant under
independent actions of spatial and microstructure flows. We showed
that this results in Doyle-Ericksen formulas identical to those
obtained from covariant energy balance for structured continua
with free microstructure manifolds.

\vskip 0.4 in
\section {\textbf{Acknowledgements}}

AY benefited from discussions with J.D. Clayton, P.M. Mariano, and
A. Ozakin.

\end{document}